\documentclass[preprint2]{aastex631}

\usepackage{amsmath,amssymb,amsfonts}
\usepackage{bm}
\usepackage{hyperref}
\usepackage[capitalise]{cleveref}
\usepackage{color}
\usepackage{soul,xcolor}
\usepackage{booktabs}

\DeclareMathOperator\vdim{dim}

\submitjournal{ApJ}

\begin{document}
\setstcolor{red}

\title{KFVM-WENO: A high-order accurate kernel-based finite volume method for compressible hydrodynamics}

\author[0000-0001-9254-7342]{Ian May}
\affiliation{Applied Mathematics, University of California Santa Cruz, \\1156 High Street, Santa Cruz, CA 95064, USA}

\author[0000-0002-8229-3040]{Dongwook Lee}
\affiliation{Applied Mathematics, University of California Santa Cruz, \\1156 High Street, Santa Cruz, CA 95064, USA}

\begin{abstract}
This paper presents a fully multidimensional kernel-based reconstruction scheme for finite volume methods applied to systems of hyperbolic conservation laws, with a particular emphasis on the compressible Euler equations. Non-oscillatory reconstruction is achieved through an adaptive order weighted essentially non-oscillatory (WENO-AO) method cast into a form suited to multidimensional reconstruction. A kernel-based approach inspired by radial basis functions (RBF) and Gaussian process (GP) modeling, which we call KFVM-WENO, is presented here. This approach allows the creation of a scheme of arbitrary order of accuracy with simply defined multidimensional stencils and substencils. Furthermore, the fully multidimensional nature of the reconstruction allows for a more straightforward extension to higher spatial dimensions and removes the need for complicated boundary conditions on intermediate quantities in modified dimension-by-dimension methods. In addition, a new simple-yet-effective set of reconstruction variables is introduced, which could be useful in existing schemes with little modification. The proposed scheme is applied to a suite of stringent and informative benchmark problems to demonstrate its efficacy and utility. A highly parallel multi-GPU implementation using Kokkos and the message passing interface (MPI) is also provided.
\end{abstract}

\section{Introduction} \label{sec:intro}
The numerical solution of systems of hyperbolic conservation laws has been a vigorous subject of research for several decades now. A notable feature of hyperbolic conservation laws is their ability to simultaneously support complicated -- but otherwise smooth -- solutions and discontinuous solutions. The nonlinear nature of the hyperbolic laws can potentially turn initially smooth flows into non-smooth flows with shocks and discontinuities. We seek to numerically solve hyperbolic conservation laws in divergence form
\begin{equation}
  \frac{\partial\bm{U}}{\partial t} + \nabla\cdot\bm{F}(\bm{U}) = \bm{S},
  \label{eq:model}
\end{equation}
where $\bm{U}$ is a vector of conserved variables, $\bm{F}$ is convex flux tensor, and $\bm{S}$ is a vector of source terms. In particular, we focus on the compressible Euler equations defined by
\begin{equation}
  \bm{U} = \begin{pmatrix} \rho \\ \rho u_i \\ E \end{pmatrix},\quad\bm{F}_j = \begin{pmatrix} \rho u_j \\ \rho u_i u_j + p\delta_{ij} \\ u_j(E + P) \end{pmatrix},
  \label{eq:euler}
\end{equation}
where the conserved quantities in $\bm{U}$ include the density, linear momenta, and total energy, respectively. This equation is closed by the adiabatic equation of state for the pressure, $p = \left(\gamma - 1\right)\left(E - \frac{1}{2}\rho u_k u_k\right)$, where $\gamma$ is the ratio of specific heats.

Numerical methods for the solution of these systems thus need to be designed around at least two competing goals: (i) accurate representation and evolution of smooth solutions and (ii) robust and physically consistent behavior in the vicinity of shocks and discontinuities. Modern computer hardware strongly favors the use of high-order accurate schemes that efficiently resolve smooth solutions via increased floating point operations and reduced loads from main memory. On the other hand, Godunov's theorem \citep{godunov1959finite} asserts that these methods must be inherently nonlinear if they are to be stable and higher than first order accurate.

Essentially non-oscillatory (ENO) schemes \citep{harten1987uniformly,harten1989eno} and subsequently weighted essentially non-oscillatory (WENO) schemes \citep{liu1994weighted,jiang1996efficient} have been one fruitful avenue for the design of methods that can simultaneously handle smooth data accurately and discontinuous data robustly. There are a plethora of different variants in WENO, which aim to improve the baseline WENO scheme of Jiang and Shu (WENO-JS) \citep{jiang1996efficient}. Some well-known examples include WENO-Z \citep{borges2008improved,castro2011high}, central-WENO \citep{levy1999central,qiu2002construction}, Hermite-WENO \citep{qiu2004hermite,qiu2005hermite}, WENO-AO \citep{balsara2016efficient}, and mapped WENO \citep{henrick2005mapped,feng2012new,feng2014improved,wang2016new,li2020modified,LI2022115004}, to list just a few. Of most relevance to the present article are the kernel-based methods of GP-WENO \citep{reyes2018new,reyes2019variable}, RBF-CWENO \citep{hesthaven:rbfcweno}, and the optimal recovery finite volume method (FVM) from \cite{sonar:OptRecFVM}.

Finite difference methods (FDMs) readily generalize to higher dimensions by applying baseline one-dimensional schemes in a dimension-by-dimension manner (e.g., see \cite{reyes2019variable}), which preserves the target one-dimensional solution accuracy in higher dimensions easily. The simplicity of achieving accuracy in higher spatial dimensions is mainly due to the fact that finite difference methods evolve pointwise quantities. In comparison, finite volume methods (FVMs) evolve volume-averaged (or cell-averaged) quantities, which complicates the design of reconstruction schemes in multiple spatial dimensions higher than second-order \citep{mccorquodale2011high,zhang2011order,balsara2009efficient,balsara2009divergence,bourgeois2022gp}. As such, numerical solutions from multidimensional finite volume methods with naive dimension-by-dimension spatial reconstruction will be limited to second-order \citep{zhang2011order,buchmuller2014improved,lee2017piecewise}.

A few key challenges arise in the development of ``multidimensional polynomial'' WENO schemes. Several practitioners have taken multivariate polynomials to design novel high-order schemes \citep{balsara2009efficient,Semplice2016}. In general, however, multidimensional interpolation/reconstruction suffers from well-known complications where not all stencils result in a solvable linear system (e.g., see Section 3.5 in \cite{reeves2022application}). Besides, the principal design choice of grid stencil requires that the full local stencil should be configured as small as possible while allowing the reconstruction of a polynomial to retain the desired degree for accuracy. The full local stencil should also be symmetric with respect to the grid to avoid having any preferred direction and allow decomposition into symmetrically placed substencils. The size and shape of the substencils also need to support the desired $k^{\rm th}$ degree multivariate polynomial reconstruction compatible with the full local stencil. Matching the full stencil and substencil sizes to the dimension of polynomial spaces while satisfying the symmetry requirements is generally impossible. All of these issues become increasingly challenging as the degree $k$ increases. Perhaps, the simplest (sub)stencil configuration would be to use rectangular stencils that allows reconstruction in a tensor product basis, although this approach comes at the cost of larger stencils than necessary to reach a given order of accuracy. Alternatively, one could reconstruct polynomials in the least squares sense by using more cells per (sub)stencil than the dimension of the corresponding space of polynomials. See more discussions in \cite{bourgeois2022gp,reeves2022application}.

Instead, one can bring one-dimensional reconstruction schemes into higher dimensions by adding transverse corrections to the baseline one-dimensional schemes for each additional direction \citep{buchHelz:dimByDim}. In two dimensions these types of schemes use cell-averaged quantities to reconstruct face-averaged quantities, which are subsequently used as input to a second reconstruction in the transverse direction to reconstruct point values on the cell faces. In three dimensions three separate reconstructions are needed (hence increasing computational expense), cascading from cell-averages to face-averages, then face-averages to line-averages, and finally line-averages to point values. Centered at cell faces, the resulting point values, referred to as Riemann states, are by design constructed to attain high-order accuracy and numerical stability in multiple dimensions. These Riemann states can be passed to any Riemann solver resulting in a pointwise-defined numerical flux at each cell face center. Similar transverse corrections can be run in reverse to generate the necessary face-averaged fluxes from point-valued fluxes. See more details in \cite{buchHelz:dimByDim}.

An alternative can be designed that only applies the transverse corrections to arrive at multiple pointwise-defined Riemann states along each transverse direction on each face and use them to subsequently evaluate the flux integrals at cell face centers through a quadrature rule \citep{zhang2011order}. As before, one uses a one-dimensional reconstruction scheme to build face-averaged quantities from the known cell-averaged quantities. In two dimensions a subsequent transverse reconstruction can be made and evaluated at several points on the face. In three dimensions the same cascade of reconstruction as described above must be done, but it needs to be done multiple times at multiple quadrature points on each 2D face to generate the appropriate pointwise-defined Riemann states.

These approaches above unfortunately suffer several essential drawbacks: (i) reconstruction of the Riemann states now requires multiple passes over the data hence requiring extra loads from memory, (ii) near the domain boundaries, the \textit{intermediate} 2D face-averaged (and edge-averaged in 3D) quantities need to be filled in accordance with the boundary conditions, which may be non-trivial, and (iii) implementations relying on distributed memory parallelism need to either perform additional communications to pass these intermediate quantities between neighboring subdomains or localize communication by utilizing much larger overlaps and repeating work adjacent to the local subdomain. In two dimensions, these drawbacks complicate the implementation mildly but are otherwise manageable. In three dimensions, however, these drawbacks become increasingly problematic. These issues are all exacerbated if adaptive mesh refinement is to be used, though it is certainly possible \citep{buchhelzdreher:amr}.

Taking a kernel-based approach to reconstruction alleviates many of the aforementioned issues. The full stencil and substencils can be chosen with \textit{size and symmetry} as the primary motives. The accuracy of reconstruction (or interpolation) on any (sub)stencil is determined by what polynomials can be \textit{reproduced} there by the chosen kernel and stencil. Kernel-based reconstruction can be seen as a generalization of the least squares approach to polynomial reconstruction, where the polynomial spaces assigned to each (sub)stencil are implicitly defined \citep{schabackWendland:KernTech,wendland:Scattered}. Fully multidimensional reconstruction is thus easy to formulate and produces pointwise-defined Riemann states that are obtained \textit{directly} from cell-average data. This removes the need for any intermediate quantities, which greatly simplifies the parallel implementation and reduces reconstruction to needing only a single pass over the data.

The present study builds around a series of recent studies exploring the use of radial basis functions (RBFs) and Gaussian Process (GP) modeling as an alternative to traditional polynomial-based approaches. The first shock-capturing finite volume GP algorithm was introduced in \cite{reyes2018new}. Called GP-WENO, this study put forward a new way of designing a versatile selectable-order property in modeling the 1D Euler and magnetohydrodynamics (MHD) equations. In addition to the newly proposed one-dimensional GP reconstruction scheme therein, this work also refactored the conventional $L_2$-based \textit{a priori} shock-handling smoothness indicators for polynomial WENO reconstruction \citep{jiang1996efficient} into an alternative that is genuinely designed for polynomial-free GP reconstruction. Similarly, the RBF-CWENO method introduced in \cite{hesthaven:rbfcweno} (and the related one dimensional scheme from \cite{bigoniHesthaven:rbfweno}) combined radial basis functions and polynomials and proposed novel smoothness indicators that treated these different components separately. The GP paradigm has also been successfully developed in a finite difference method (FDM). Interested readers can refer to the work by  \cite{reyes2019variable}, where a full three-dimensional conservative GP-WENO finite difference scheme solves a stringent set of benchmark problems in the system of the Euler equations with improved solution accuracy and computational efficiency.

It is worth mentioning that the present study will be limited to the use of uniform Cartesian grids. Obviously, much simpler and cheaper finite difference methods (FDMs) could be used on these structured grids. We acknowledge that our numerical examples presented herein can be solved with lighter computational workload with FDMs.
Regardless, we emphasize that the design of multidimensional finite volume methods remains a necessary task as they are much easier to incorporate into solvers utilizing adaptive mesh refinement (AMR) and into solvers utilizing unstructured meshes. Both AMR and unstructured meshes lie beyond the scope of the present study, but the presented scheme is designed with possible extensions to these cases in mind. The few items that truly rely on the underlying uniform mesh will be called out. Naturally, demonstrating the efficacy of the proposed scheme in the simpler case of a uniform mesh is a necessary step towards constructing more complex solvers in the future.

The aforementioned schemes all employ \textit{a priori} nonlinear limiting which aims to suppress and avoid numerical oscillations before they can form. Unlike these \textit{a priori} GP-WENO methods, a two-dimensional \textit{a-posteriori} GP-MOOD (Multidimensional Optimal Order Detection) scheme in finite volume has been studied and reported in \cite{bourgeois2022gp}. GP-MOOD adopts a set of high-order GP reconstruction methods (e.g., $3^{\rm rd}$, $5^{\rm th}$, and $7^{\rm th}$) in place of the polynomial reconstructions in the original polynomial-MOOD approach \citep{clain2011high,diot2012improved,diot2013multidimensional,diot2012methode}. The resulting method is a strong positivity-preserving GP-based solver for shock-dominant compressible flows. This new GP-MOOD further improves upon the ingredients of the original MOOD limiting in the polynomial MOOD methods \citep{clain2011high,diot2012improved,diot2013multidimensional,diot2012methode} by introducing the new Compressibility-Shock Detection (CSD) switch that controls a good balance between numerical accuracy and diffusivity. Besides, this new CSD switch can improve Discontinuous Galerkin (DG) methods that adopt the polynomial MOOD shock detection as subcell-based limiting to switch from DG to an alternative shock-stable method (e.g., FVM) at shock cells
\citep{dumbser2014posteriori}.

There have been other related GP applications beyond its role as a hyperbolic solver. In the recent work by \cite{reeves2022application}, the GP-WENO method was extended to a 3rd-order prolongation algorithm in finite volume adaptive mesh refinement (AMR) simulations using AMReX \citep{zhang2019amrex}. As another non-fluid dynamics application, GP interpolation was shown as an improved mathematical tool for upsampling optical characters from coarse resolutions to fine resolutions \citep{reeves2020gaussian}.

This paper introduces the KFVM-WENO (kernel-based finite volume method with WENO) scheme and proceeds as follows. Section \ref{sec:recon} develops the kernel-based reconstruction scheme that is central to this work. Section \ref{sec:weno} incorporates WENO into the reconstruction scheme to handle discontinuous data, deliberates on the choice of variables for reconstruction, and closes with an adaptation of the KXRCF indicator \citep{krivodonova:KXRCF} to flag when the nonlinear WENO limiting is required. Section \ref{sec:pos_pres} documents the incorporation of the self-adjusting positivity-preserving limiter detailed in \cite{balsara:Astrojets}. Finally, the method in full is reviewed in Section \ref{sec:method}, followed by a suite of benchmark problems in Section \ref{sec:results}. 

\section{Kernel-based reconstruction} \label{sec:recon}
The first major component of a finite volume method for the solution of systems of conservation laws is the reconstruction of the state within each cell, which can then be evaluated to find high-order accurate pointwise-valued Riemann states along the boundary of the cell (or cell interface). This reconstruction, $\widetilde{U}(x)$, should simultaneously provide an accurate approximation of the true state and remain well-behaved in the presence of shocks or other discontinuities. In the present study, we pose this reconstruction problem as optimal recovery in a reproducing kernel Hilbert space (RKHS) \citep{hesthaven:rbfcweno,schabackWendland:KernTech,sonar:OptRecFVM,rasmussen2005}. 
This approach yields a form that is mostly dimension-independent and allows great flexibility in the choice of stencils. We note that related but distinct formulations have been examined previously \citep{reyes:gpWenoFvm,guoJung:RbfWENO,liu:wlsENO,aboiyar:RBFADER}.

By convention, we refer to this RKHS approximation as \textit{reconstruction} (or generalized interpolation) when designed to match cell averages, and note that other more general input data could also be supplied. In this section, a kernel-based method for reconstruction is presented. The reconstruction presented here is \textit{linear} (i.e., lacks nonlinear limiting) with respect the local stencil input data only linearly without limiting, and is hence inappropriate for use near shocks. This is resolved in Section \ref{sec:weno}.

\subsection{Asymmetric reconstruction} \label{sec:recon:asym}
Consider a set of finite volumes $\Omega_h\subset\mathbb{R}^d$ for $h=1,\ldots,N$, each with volume $||\Omega_h||$ measured in the standard Eucliean $L_2$ norm. Let us denote the cell averaging functional for $\Omega_h$ with respect to $\bm{x} \in \Omega_h$ by
\begin{equation}
  \lambda_h^{(\bm{x})} = \frac{1}{||\Omega_h||}\int\limits_{\Omega_h} \bm{\cdot}~d\bm{x}.
  \label{eq:lambda_avg}
\end{equation}
The given data for reconstruction are cell averages of $f(\bm{x})$, which we gather into the vector $\bm{g}$ whose entries are ${\rm g}_h = \lambda_h^{(\bm{x})}f(\bm{x})$. 
Asymmetric kernel \textit{reconstruction} seeks an approximant $\widetilde{f}(\bm{x})$ of the form
\begin{equation}
  \label{eq:avExpan}
  \widetilde{f}(\bm{x}) = \sum\limits_{l=1}^N a_l K(\bm{x},\bm{x}_l) + \sum\limits_{|\bm{\alpha}_{(v)}|\leq p} b_v \bm{x}^{\bm{\alpha}_{(v)}},
\end{equation}
which is a highly accurate approximation to the underlying function $f(\bm{x})$. More specifically, we aim to find $\widetilde{f}(\bm{x}_*) \approx {f}(\bm{x}_*)$ that represents pointwise Riemann states at $\bm{x}_*$.

Here, $K(\bm{x},\bm{y}):\mathbb{R}^d\times\mathbb{R}^d\rightarrow\mathbb{R}$ is a symmetric positive definite kernel function \citep{schabackWendland:KernTech,fasshauer2015kernel}, and $\bm{x}_l$ denotes the center of the $l^{\rm th}$ cell. In this work we take the kernel function to be the squared exponential
\begin{equation}
  \label{eq:kernSE}
  K(\bm{x},\bm{y}) = e^{-\frac{||\bm{x} - \bm{y}||^2}{2\ell^2}},
\end{equation}
where the hyperparameter $\ell$ is a length scale.

The second summation in \cref{eq:avExpan} augments the kernel expansion in the first summation with an additional polynomial expressed in the monomial basis. The monomials $\bm{x}^{\bm{\alpha}_{(v)}}$ span the space $\mathbb{P}_p^d$ of polynomials with maximum total degree $p$ in $d$ spatial dimensions. Each $\bm{\alpha}_{(v)}$ is a multi-index (a $d$-tuple $\bm{\alpha}\in\mathbb{N}^d_0$ applying componentwise as a power and having $|\bm{\alpha}|=\alpha_1 + \alpha_2 + \ldots$) where the subscript $(v)$ indexes the $D=\vdim\left(\mathbb{P}_p^d\right)$ separate monomials. Here and throughout, parenthetical subscripts on boldface symbols (e.g., ${(v)}$ on $\bm{\alpha}_{(v)}$) are meant to label one such vector or multi-index from a larger collection, while bare subscripts on plain symbols (e.g., $h$ on $g_h$) indicate individual entries of a vector, matrix, or multi-index.

The polynomial term in \cref{eq:avExpan} is not strictly required when using kernels such as the squared exponential, though it would be necessary for well-posedness of the reconstruction problem when using conditionally positive definite kernels such as the thin-plate splines \citep{schabackWendland:KernTech,flyer:RbfWithPolys}. However, the inclusion of this polynomial term is still useful as it allows for the asymptotic accuracy of the reconstruction process to be maintained independent of the hyperparameter $\ell$ \citep{flyer:RbfWithPolys}. This allows $\ell$ to be fixed at moderate lengths to avoid problems with extremely ill-conditioned linear systems. The inclusion of the polynomial term forgoes any need for temporary elevations of precision as present in \cite{reyes2018new} or the use of elaborate stable inversion schemes \citep{fasshauer2012stable,fornberg:rbfpade,wright:vecvalrational}.

Two important observations are in order. First, the polynomial term does not restrict stencil selection in the same way that is seen in pure high-order piecewise polynomial schemes. The only requirement is that the included space of polynomials be unisolvent over the stencil, and in practice this means picking a stencil first and subsequently choosing the maximum polynomial degree $p$ to maintain unisolvence (see \cite{flyer:RbfWithPolys} for a more fundamental discussion of this approach). Second, the form of the approximant $\widetilde{f}(\bm{x})$ in \cref{eq:avExpan} is very generic and places little to no restriction on the stencil layout or shape of cells. These properties provide the choice of KFVM-WENO's multidimensional stencil with great flexibility (see more discussions in Section \ref{sec:weno:stencils}), relieving us from the stencil-related issues described in Section \ref{sec:intro}.

Enforcing that \cref{eq:avExpan} matches the given cell-averaged data 
$\bm{g} = (g_1, \dots, g_h, \dots, g_N)^T$, 
\begin{equation}
\lambda_h^{(\bm{x})} \widetilde{f}(\bm{x}) = {\rm g}_h,\quad h=1,\ldots,N,
\end{equation}
requires that the coefficient vectors $\bm{a}$ and $\bm{b}$ satisfy
\begin{equation}
  \bm{Q}\bm{a} + \bm{P}\bm{b} = \bm{g} \label{eq:kernAppxCond}.
\end{equation}
The entries of the $N \times N$ kernel matrix $\bm{Q}$ correspond to the integrals of the kernel function $K$ anchored at each $\bm{x}_l$, $l = 1, \dots, N$, over each cell $\Omega_h$, i.e.,
\begin{equation}
{\rm Q}_{hl} 
= \lambda_h^{(\bm{x})}K(\bm{x},\bm{x}_l) \label{eq:kernelMat} 
=\frac{1}{||\Omega_h||} \int\limits_{\Omega_h}K(\bm{x},\bm{x}_l) d\bm{x}.
\end{equation}
Similarly, the entries of the $N \times D$ matrix $\bm{P}$ are the integrals of the monomials $\bm{x}^{\bm{\alpha}_{(v)}}$ over each cell $\Omega_h$, i.e.,
\begin{equation}
{\rm P}_{hv} 
= \lambda_h^{(\bm{x})}\bm{x}^{\bm{\alpha}_{(v)}}
=\frac{1}{||\Omega_h||} \int\limits_{\Omega_h} \bm{x}^{\bm{\alpha}_{(v)}} d\bm{x}.
\end{equation}

At present the linear system \cref{eq:kernAppxCond} is under-determined since there are $N+D$ unknowns in $\bm{a}$ and $\bm{b}$ but only $N$ equations. To resolve this issue, we follow the same approach used in the RBF literature (e.g., see \cite{bigoni2017adaptive}) and further require that \cref{eq:avExpan} exactly reproduces polynomials spanned by the monomials present in the second summation. More specifically, this means that the set of coefficients $\bm{a} \in \mathbb{R}^N$ must lie orthogonal to the polynomial space, and hence must satisfy $\bm{P}^T\bm{a} = \bm{0}$. This forces the first summation in \cref{eq:avExpan} to only fit data that lies outside the span of the monomial terms. Putting things together, the coefficients $\bm{a}$ and $\bm{b}$ are found by solving the $(N+D)\times (N+D)$  block linear system
\begin{equation}
  \label{eq:recBlockSys}
  \begin{bmatrix} \bm{Q} & \bm{P} \\ \bm{P}^T & \bm{0} \end{bmatrix} \begin{pmatrix} \bm{a} \\ \bm{b} \end{pmatrix} = \begin{pmatrix} \bm{g} \\ \bm{0} \end{pmatrix}.
\end{equation}

Returning to \cref{eq:avExpan} 
assuming $\bm{a}$ and $\bm{b}$ are now known,
the approximant $\widetilde{f} \approx {f}$ can be evaluated at a point $\bm{x}_*$ via
\begin{align}
  \widetilde{f}\left(\bm{x}_*\right)
  &= \bm{T}^T\bm{a} + \bm{S}^T\bm{b} \label{eq:kernAppxEval} \\
  &=  \begin{pmatrix} \bm{T}^T | \; \bm{S}^T \end{pmatrix} \begin{pmatrix} \bm{a} \\ \bm{b} \end{pmatrix}
  \label{eq:kernAppxEval2}  \\
  &= \begin{pmatrix} \bm{T}^T | \; \bm{S}^T \end{pmatrix}  \begin{bmatrix} \bm{Q} & \bm{P} \\ \bm{P}^T & \bm{0} \end{bmatrix}^{-1} \begin{pmatrix} \bm{g} \\ \bm{0} \end{pmatrix},\label{eq:kernAppxEvalBlock}
\end{align}
where the entries of $\bm{T}$ and $\bm{S}$ are respectively given by
\begin{align}
{\rm T}_l &= K(\bm{x}_*,\bm{x}_l) \label{eq:sampleVec}, \;\; 1 \le l \le N,\\
{\rm S}_v &= \bm{x}_*^{\bm{\alpha}_{(v)}}, \;\; 1 \le v \le D.
\end{align}
However, building and solving the block system in \cref{eq:recBlockSys} every time reconstruction is needed would be very computationally expensive. Fortunately, as hinted by \cref{eq:kernAppxEvalBlock}, most of the work being done in this process is independent of the data $\bm{g}$ that varies temporally and spatially. Considering the terms other than the data vector $\bm{g}$ yields that a reconstruction vector $\bm{r}$ could be instead computed as the solution of
\begin{equation}
  \label{eq:recVec1}
  \begin{bmatrix} \bm{Q}^T & \bm{P} \\ \bm{P}^T & \bm{0} \end{bmatrix} \begin{pmatrix} \bm{r} \\ \bm{w} \end{pmatrix} = \begin{pmatrix} \bm{T} \\ \bm{S} \end{pmatrix},
\end{equation}
which gives
\begin{align}
   \begin{pmatrix} \bm{r} \\ \bm{w} \end{pmatrix}^T&=
  \begin{pmatrix} \bm{T} \\ \bm{S} \end{pmatrix}^T 
  \begin{bmatrix} \bm{Q}^T & \bm{P} \\ \bm{P}^T & \bm{0} \end{bmatrix}^{-T}   
  \label{eq:recVec2}\\ &=
  \begin{pmatrix} \bm{T}^T | \;  \bm{S}^T \end{pmatrix}
  \begin{bmatrix} \bm{Q} & \bm{P} \\ \bm{P}^T & \bm{0} \end{bmatrix}^{-1}.
    \label{eq:recVec3}
\end{align}
This reduces the evaluation of the approximant in \cref{eq:kernAppxEvalBlock} to a
single dot product
\begin{equation}
  \widetilde{f}\left(\bm{x}_*\right)
  = \begin{pmatrix} \bm{r} \\ \bm{w} \end{pmatrix}^T  
  \begin{pmatrix} \bm{g} \\ \bm{0} \end{pmatrix} 
  = \bm{r}^T\bm{g} 
  = \bm{r} \cdot \bm{g}.
  \label{eq:reconVec}
\end{equation}
Note that the weights $\bm{w}$ associated with the monomial terms are not explicitly needed for anything. This construction can be viewed as an optimization problem in the associated RKHS wherein $\bm{w}$ are the Lagrange multipliers from the polynomial exactness constraints. Since $\bm{w}$ is not directly used, the cost of producing one such point value depends only on the stencil size $N$ and is independent of the maximal polynomial degree $p$.

Each evaluation point $\bm{x}_*$ will need its own reconstruction vector. The introduction of WENO in Section \ref{sec:weno} will further require different reconstruction vectors for each substencil. Once a mesh configuration is determined, the reconstruction vectors $\bm{r}$ can be easily precomputed in an initialization step for our KFVM-WENO solver as a whole. This means that reconstruction to a point value only requires a single dot product (see \cref{eq:reconVec}) between vectors whose size is determined by the stencil. For the uniform grids considered here we only need to compute one set of reconstruction vectors for a generic stencil and its substencils, detailed in Section \ref{sec:weno}. For unstructured grids a set of reconstruction vectors would need to be cached for each cell in the mesh as no two cells are likely to have identical stencils. Fortunately, the calculation of the various reconstruction vectors is trivially parallelizable and remains something that needs to be done \textit{only once} during initialization.

\subsection{Ties to Gaussian Processes} \label{sec:recon:gp}
We note that this reconstruction method bears similarity to evaluations of the \textit{updated posterior predictive mean} of a Gaussian process (GP) in the zero-noise limit \citep{reyes:gpWenoFvm,reyes2019variable,bourgeois2022gp,fasshauer:KernApprox,rasmussen2005,bishop2007pattern}. Note that in general 
\begin{equation}
{\rm Q}_{hl} =
\lambda_h^{(\bm{x})}K(\bm{x},\bm{x}_l) \neq \lambda_l^{(\bm{x})}K(\bm{x},\bm{x}_h) = {\rm Q}_{lh} , 
\end{equation}
hence the matrix $\bm{Q}$ in \cref{eq:kernelMat} is not symmetric. However, cell averaging functionals $\lambda_h^{(\bm{x})}$ could also be incorporated into the kernel expansion in \cref{eq:avExpan}, thus symmetrizing the arising system. In doing so $\bm{Q}$ can be interpreted as a covariance matrix for a Gaussian process as was done to great success in \cite{reyes:gpWenoFvm,reyes2019variable,bourgeois2022gp}.
 
We make another consequential remark on GP reconstruction. Viewing the reconstruction procedure through the lens of Gaussian processes could provide insight into the optimization of the length scale $\ell$, and possibly into uncertainty quantification via an updated posterior covariance kernel \citep{rasmussen2005,bishop2007pattern}. As length scale optimization and uncertainty quantification lie beyond the scope of this work, the remainder of this paper will only concern the \textit{deterministic} interpretation of kernel reconstruction presented in the preceding section.

\section{Nonlinear reconstruction} \label{sec:weno}
The reconstruction scheme discussed so far is \textit{linear} with respect to the input data, which will unavoidably generate unwanted oscillations in the vicinity of shocks and discontinuities. In this section, a weighted essentially non-oscillatory (WENO) scheme is developed for these kernel-based reconstructions to control such oscillations. The great flexibility of kernel-based reconstruction allows many different stencil and substencil configurations.

\subsection{Stencils and substencils} \label{sec:weno:stencils}
\begin{figure}
  \centering
  \includegraphics[width=0.9\linewidth]{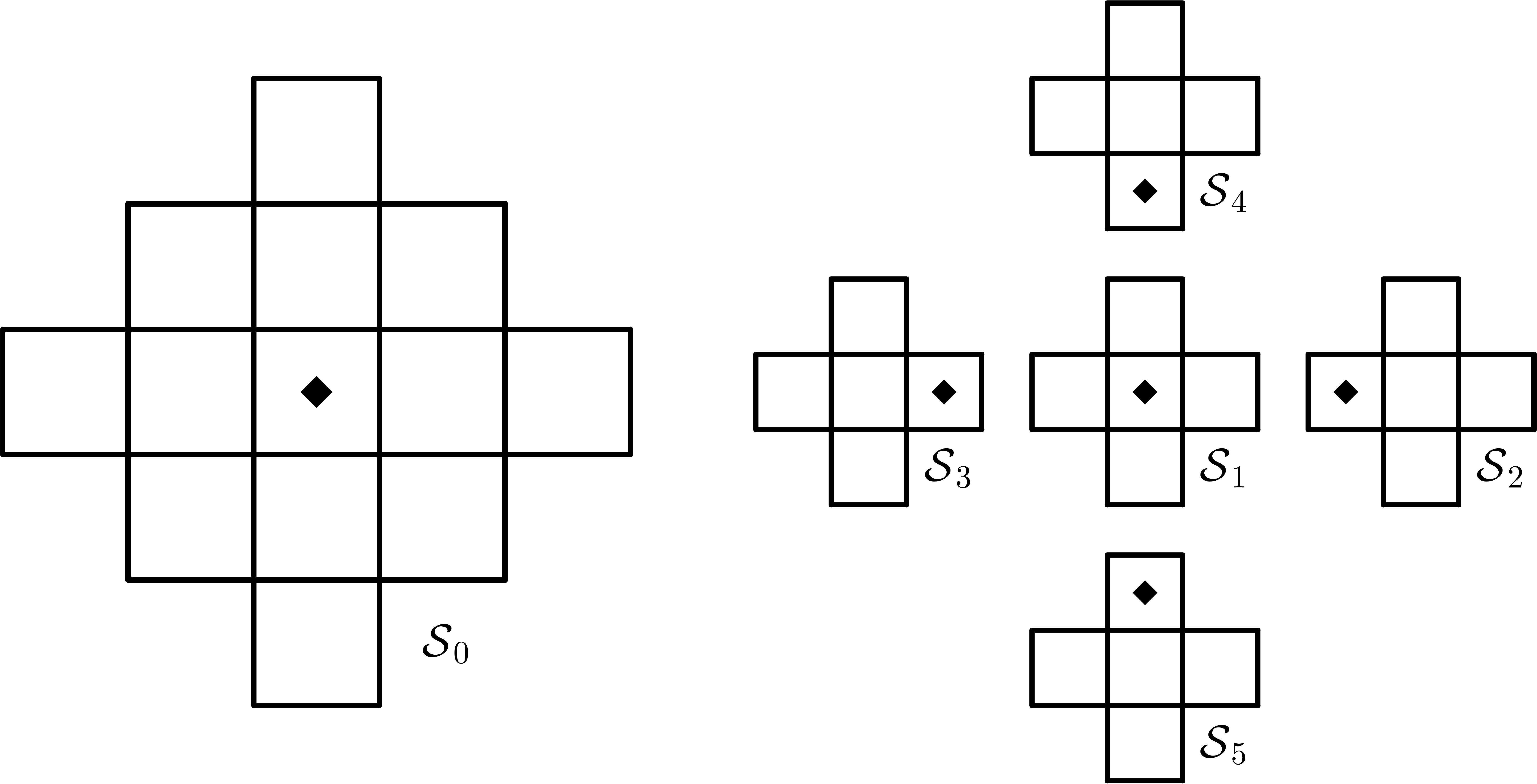}
  \caption{The full radius-2 stencil $S_0$ is shown on the left with its five substencils, $S_q, q=1, \dots, 5$, on the right. The cell with the diamond in each (sub) 
stencil indicates the central cell containing $\bm{x}_*$ where reconstruction is performed. The full stencil $S_0$ has 13 cells, while each of the substencils has five.}
  \label{fig:rad2}
\end{figure}

\begin{figure}
  \centering
  \includegraphics[width=0.9\linewidth]{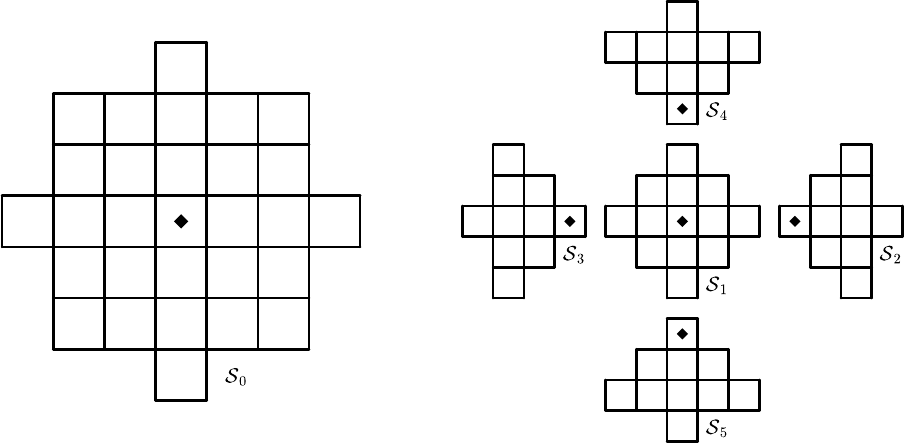}
  \caption{The full radius-3 stencil $\mathcal{S}_0$ is shown on the left with its five substencils, $\mathcal{S}_q, q=1, \dots, 5$, on the right. The cell with the diamond in each (sub)stencil indicates the central cell containing $\bm{x}_*$ where reconstruction is being performed. The full stencil $\mathcal{S}_0$ has 29 cells, the central substencil $\mathcal{S}_1$ has 13 cells, and each of the remaining biased substencils, $\mathcal{S}_2, \dots, \mathcal{S}_5$, has 10.
  }
  \label{fig:rad3}
\end{figure}

We consider ``spherical'' stencils of radius $R=2,3$, e.g., see \cref{fig:rad2,fig:rad3}. These stencils are unisolvent (or uniquely solvable) over monomials of total 
degree $3$ and $5$ respectively and hence yield reconstructions accurate at least to order $\mathcal{O}(\Delta^{2R})$, where $\Delta$ is the grid spacing. 
Careful selection of the hyperparameter $\ell$ can yield higher accuracy, tending towards $\mathcal{O}(\Delta^{2R+1})$, which we discuss further in Section \ref{sec:method:precomp}. Examples of two dimensional stencils can be seen in \cref{fig:rad2} and \cref{fig:rad3} for $R=2$ and $R=3$, respectively.

Let $\mathcal{S}_0 = \lbrace (i,j)\in\mathbb{Z}^2 : (i^2+j^2) \leq R^2 \rbrace$ denote the full circular stencil of radius $R$, indexed relative to the central cell 
where $(i,j)=(0,0)$ and reconstruction is being performed. This full stencil $S_0$ is then broken into $N_S=5$ substencils $S_q$, including one central substencil 
of smaller radius, $\mathcal{S}_1$, and four other biased substencils, $\mathcal{S}_2, \dots, \mathcal{S}_5$, given as:
\begin{align*}
  \mathcal{S}_1 &= \lbrace (i,j)\in\mathcal{S}_0 : (i^2+j^2) \leq (R-1)^2 \rbrace, \\
  \mathcal{S}_2 &= \lbrace (i,j)\in\mathcal{S}_0 : |j| \leq i \rbrace, \\
  \mathcal{S}_3 &= \lbrace (i,j)\in\mathcal{S}_0 : |j| \leq -i \rbrace, \\
  \mathcal{S}_4 &= \lbrace (i,j)\in\mathcal{S}_0 : |i| \leq j\rbrace, \\
  \mathcal{S}_5 &= \lbrace (i,j)\in\mathcal{S}_0 : |i| \leq -j\rbrace.
\end{align*}
Stencils in three dimensions are formed in precisely the same way. First the full stencil is set as $\mathcal{S}_0 = \lbrace (i,j,k)\in\mathbb{Z}^3 : (i^2+j^2+k^2) \leq 
R^2 \rbrace$, again indexed relative to the central cell where $(i,j,k)=(0,0,0)$. This stencil is then split into $N_S=7$ substencils, including one central substencil 
of smaller radius, $\mathcal{S}_1$, and six other biased stencils, $\mathcal{S}_2, \dots, \mathcal{S}_7$, given as:
\begin{align*}
  \mathcal{S}_1 &= \lbrace (i,j,k)\in\mathcal{S}_0 : (i^2+j^2+k^2) \leq (R-1)^2 \rbrace, \\
  \mathcal{S}_2 &= \lbrace (i,j,k)\in\mathcal{S}_0 : |j| \leq i,~|k| \leq i \rbrace, \\
  \mathcal{S}_3 &= \lbrace (i,j,k)\in\mathcal{S}_0 : |j| \leq -i,~|k| \leq -i \rbrace, \\
  \mathcal{S}_4 &= \lbrace (i,j,k)\in\mathcal{S}_0 : |i| \leq j,~|k| \leq j \rbrace, \\
  \mathcal{S}_5 &= \lbrace (i,j,k)\in\mathcal{S}_0 : |i| \leq -j,~|k| \leq -j \rbrace, \\
  \mathcal{S}_6 &= \lbrace (i,j,k)\in\mathcal{S}_0 : |i| \leq k,~|j| \leq k \rbrace, \\
  \mathcal{S}_7 &= \lbrace (i,j,k)\in\mathcal{S}_0 : |i| \leq -k,~|j| \leq -k \rbrace.
\end{align*}

For a given point $\bm{x}_*$ in the cell where point values are sought, \cref{eq:kernAppxCond} through \cref{eq:reconVec} furnish a different reconstruction vector for each (sub)stencil. We denote the $q^{\rm th}$ such reconstruction vector by $\bm{r}_{(q)}$ on $S_q$ for each $q=0, \dots, N_S$. 

Note that the (sub)stencils, as written, consist of collections of unit cells (i.e., cell size being unity) rather than cells matching the grid scale. The integrals required for filling the matrix $\bm{Q}$ in \cref{eq:kernelMat} and the sample vector $\bm{T}$ in \cref{eq:sampleVec} can easily be calculated over these unit cells. This change of variables only requires that the length scale $\ell$ be replaced by $\ell/\Delta$ inside the kernel and naturally cancels out the pre-factor that scales by the cell volume and decouples the calculation of the reconstruction vectors $\bm{r}_{(q)}$ from the description of the grid. Additionally, this reformulation will be exploited to generate grid-independent smoothness indicators in Section \ref{sec:weno:nonlin}.

Furthermore, these particular choices of substencils are the only remaining item that is specific to uniform grids, though there would be no great difficulty in selecting substencils on unstructured grids.

\subsection{Adaptive order WENO} \label{sec:weno:ao}
A traditional WENO method leverages only the substencils (i.e., excluding the full stencil $\mathcal{S}_0$) to recover the action of the full stencil $\mathcal{S}_0$ \textit{indirectly} through cleverly chosen nonlinear weights, $\omega_q$, that tend towards appropriately chosen linear weights, $\gamma_q$, when the contained data is smooth. Mathematically speaking, recovering the action of the full stencil is only possible when the $0^{\rm th}$ reconstruction vector $\bm{r}_{(0)}$ on $\mathcal{S}_0$ lies in the span of the remaining reconstruction vectors, i.e., $\bm{r}_{(0)} \in \text{span}\{\bm{r}_{(q)}: q=1, \dots, N_S\}$, assuming the reconstruction vectors of the substencils, $\bm{r}_{(q)}, q=1, \dots, N_S$, have been appropriately extended by zero-padding to match the dimension of $\bm{r}_{(0)}$. The exact linear weights (if found) would then satisfy $\bm{r}_{(0)} = \sum\limits_{q=1}^{N_S}\gamma_q \bm{r}_{(q)}$.
Unfortunately, however, the full stencil reconstruction vector $\bm{r}_{(0)}$ does not lie in the span of the remaining reconstruction vectors. 

Instead, inexact linear weights could be found by solving a least-squares problem (e.g., see \cite{reyes2018new} for a 1D GP-WENO in FVM and \cite{reyes:gpWenoFvm} for a multidimensional GP-WENO in FDM). Unfortunately, though, we have observed a degradation in accuracy when applying this approach to multidimensional FVM reconstruction; hence the least-squares approach is not appropriate for the current study.

Alternatively, the adaptive order WENO method (WENO-AO) employs the full stencil $\mathcal{S}_0$ \textit{directly} and selects weights solely to provide stability \citep{balsara:wenoAO}. Given a vector of cell averages $\bm{g}_{(k)}$ of $f(\bm{x})$ supported on the $k^{\rm th}$ (sub)stencil, KFVM-WENO forms the reconstruction (see also Appendix A in \cite{reyes2019variable}),
\begin{equation}
  \widetilde{f}(\bm{x}_*) = \frac{\omega_0}{\gamma_0}\bm{r}_{(0)}^T\bm{g}_{(0)} + \sum\limits_{q=1}^{N_S} \left(\omega_q - \omega_0\frac{\gamma_q}{\gamma_0}\right)\bm{r}_{(q)}^T\bm{g}_{(q)}.
  \label{eq:wenoAOrec}
\end{equation}
Note that if $\omega_q\rightarrow\gamma_q$ for all $q$ in the absence of discontinuities on the highest order full stencil $\mathcal{S}_0$, the coefficient on the first term tends to one, and the term in parentheses tends to zero. Hence, when the nonlinear weights $\omega_q$ approach the linear weights $\gamma_q$, the highest accuracy reconstruction $\bm{r}_{(0)}^T\bm{g}_{(0)}$ is selected, hence full accuracy can be maintained independently of how $\gamma_q$ are chosen. Alternatively, if some substencils contain a discontinuity then those weights will fall to zero, thereby minimizing the contribution from those substencils. In this case, the full stencil will obviously also contain the discontinuity so $\omega_0\rightarrow 0$, leaving only the substencils with smooth data to participate in the reconstruction. 

In two space dimensions, we follow \cite{balsara:wenoAO} to fix $\gamma_{hi} = \gamma_{lo} = 0.8$ and set the linear weights as $\gamma_0 = \gamma_{hi}$ for $\mathcal{S}_0$, $\gamma_1 = (1-\gamma_{hi})\gamma_{lo}$ for $S_1$, and $\gamma_q = (1-\gamma_{hi})(1-\gamma_{lo})/4$ for the remaining biased substencils $\mathcal{S}_q$ with $q=2, \dots, 5$. In three dimensions, the linear weights are set in the same way apart from the biased weights that become $\gamma_q = (1-\gamma_{hi})(1-\gamma_{lo})/6$ with $q=2, \dots, 7$. We provide the explicit form of the nonlinear weights $\omega_q$ in the next section.

\subsection{Smoothness indicators and weights} \label{sec:weno:nonlin}
The nonlinear weights, including unnormalized $\widetilde{\omega_q}$ and the associated normalized $\omega_q$, and the smoothness indicators $\beta_q$ within them, must be specified in an appropriate way for the reconstruction in \cref{eq:wenoAOrec} to simultaneously provide accurate reconstructions for smooth data and limited reconstructions for non-smooth data. In this section, we introduce \textit{multidimensional} smoothness indicators $\beta_q$ on each $S_q$ for our KFVM-WENO, which measure the relative smoothness of the reconstructed function $\widetilde{f}_q$ in \cref{eq:avExpan} over each $S_q$ in a \textit{multidimensional} way constructively.

We design these indicators $\beta_q$ so that they approximate the square of a scaled $H_p$ semi-norm $|\widetilde{f}_q|_{H_p}^2$, i.e., $\beta_q \approx |\widetilde{f}_q|_{H_p}^2$, where we set
\begin{align}
  |\widetilde{f}_q|_{H_p}^2 =& \sum\limits_{0 < |\bm{\alpha}|\leq p} \Delta^{2|\bm{\alpha}|-1} \int\limits_{\Omega} \left(\frac{\partial^{|\bm{\alpha}|}\widetilde{f}_q}{\partial \bm{x}^{\bm{\alpha}}} \right)^2 d\bm{x} 
  \label{eq:scaledH2_a} \\ 
  =& \sum\limits_{0 < |\bm{\alpha}|\leq p} \Delta^{2|\bm{\alpha}|} \int\limits_{\hat{\Omega}} \left( \frac{\partial^{|\bm{\alpha}|}\widetilde{f}_q}{\partial \hat{\bm{x}}^{\bm{\alpha}}} \right)^2 d\hat{\bm{x}}.
  \label{eq:scaledH2}
\end{align}
Here, $\Omega$ is the central cell where reconstruction is being performed and $\bm{\alpha}$ is a multi-index. The factor $\Delta^{2|\bm{\alpha}|-1}$ is present in the first formula \cref{eq:scaledH2_a} to make the quantity in the parentheses independent of the grid scale mimicking the use of undivided differences in standard polynomial WENO methods, see for example, \cite{jiang1996efficient}. The second line \cref{eq:scaledH2} makes the change of variables $\bm{x}=\hat{\bm{x}}\Delta$ (whereby ${\hat{\Omega}}$ is introduced accordingly), which absorbs one factor of $\Delta$ and places this in a formulation similar to the discussion regarding reconstruction vectors in Section \ref{sec:weno:stencils}.

The scaled $H_p$ semi-norm in the second part \cref{eq:scaledH2} cannot be integrated in closed form. To define $\beta_q$ such that $\beta_q \approx |\widetilde{f}_q|_{H^2}^2$, we use a simple midpoint quadrature rule
\begin{equation}
  \beta_q = \sum\limits_{0 < |\bm{\alpha}|\leq p} \Delta^{2|\bm{\alpha}|} \left( \left. \frac{\partial^{|\bm{\alpha}|}\widetilde{f}_q}{\partial \hat{\bm{x}}^{\bm{\alpha}}} \right|_{\bm{0}} \right)^2,
  \label{eq:smoothIndic}
\end{equation}
where the coordinate system has been shifted to place the center of the cell  where reconstruction is occurring to the origin as denoted by ${\bm{0}}$ in \cref{eq:smoothIndic}.

The evaluation of the partial derivatives of $\widetilde{f}_q$ in \cref{eq:smoothIndic} from given cell-average data $\bm{g}_{(q)}$ on the $q^{\rm th}$ (sub)stencil proceeds in much the same way in the reconstruction of Riemann states along the cell interface as previously done in \cref{eq:kernAppxEval,eq:kernAppxEval2,eq:kernAppxEvalBlock}, except that we now need a \textit{differential} version to compute those partial derivatives instead of the previous version in \cref{eq:sampleVec}. This can be accomplished by replacing the integral version of the sample vector $\bm{T}$ in \cref{eq:sampleVec} by a new differential version of the sample vector
\begin{equation}
  \bm{D}_l = \left.\frac{\partial^{|\bm{\alpha}|}K(\hat{\bm{x}},\hat{\bm{x}}_l)}{\partial\hat{\bm{x}}^{\bm{\alpha}}}\right|_{\bm{0}},
\end{equation}
and similarly replacing $\bm{S}$ by derivatives of the monomial terms.
In addition, we follow a similar process in \cref{eq:recVec1,eq:recVec2,eq:recVec3} to obtain a separate reconstruction vector $\bm{r}$ required for each partial derivative $\partial^{|\bm{\alpha}|}/\partial\bm{x}^{\bm{\alpha}}$. 

Increasing the maximum order of differentiation, $p$ in \cref{eq:scaledH2}, improves the detection of discontinuous data and gives a more diffusive scheme overall. However, the calculation of these smoothness indicators is the most expensive part of the reconstruction process, and hence $p$ should be chosen as small as possible while still avoiding the appearance of unwanted oscillations. Here we set $p = R$ yielding 5 (or 9) derivatives in two dimensions for radius 2 (or 3) stencils and 9 (or 19) derivatives in three dimensions for radius 2 (or 3) stencils. We have observed that these degrees of derivatives are fully sufficient in suppressing oscillations, and higher degrees induce higher computational costs for little to no benefit.

Finally, the nonlinear weights $\omega_q$ are generated from the smoothness indicators $\beta_q$ in the standard WENO-JS form \citep{jiang1996efficient}. The unnormalized weights for each stencil are calculated as
\begin{equation}
  \widetilde{\omega_q} = \frac{\gamma_q}{\beta_q^2 + \epsilon},
  \label{eq:nlinWts}
\end{equation}
where $\epsilon$ is included to avoid division by zero. Finally, the unnormalized nonlinear weights $\widetilde{\omega_q}$ are normalized to obtain the so-called normalized weights as
\begin{equation}
  \omega_q = \frac{\widetilde{\omega_q}}{\sum\limits_q \widetilde{\omega_q}}.
\end{equation}

There are several ways to set a value for $\epsilon$. First, $\epsilon$ can be set to a fixed small number independent of the grid scale as is done in many WENO schemes (see e.g., \cite{jiang1996efficient}). In contrast, numerous follow-up studies have shown that $\epsilon$ must vary with the grid scale to maintain high order accuracy for many schemes, particularly in the resolution of critical points (see e.g., \cite{henrick2005mapped,craveroCWENO:2018}). Here we take the former route and set $\epsilon = 10^{-40}$ which provides excellent robustness in the vicinity of very strong shocks. Meanwhile, we also have tested a grid dependent scaling satisfying $\epsilon \sim \Delta^R$ and observed improved accuracy in some cases, which did not contribute to further enhancement in general. Hence a constant $\epsilon$ a chosen in the current study.
%
Moreover, in Section \ref{sec:weno:kxrcf}, we introduce indicators that allow nonlinear limited WENO reconstruction to be bypassed entirely and hence the simpler and more robust choice of our single fixed $\epsilon$ will not degrade the accuracy of the method as a whole.

\subsection{Variables for reconstruction} \label{sec:weno:vars}
The multidimensional WENO reconstruction procedure detailed in the prior subsections is only defined for scalar fields, thus reconstructing the Riemann states requires that the cell average states be separated into scalar fields. It is well known that performing \textit{limited} reconstruction directly on the conservative variables leads to oscillatory results in the vicinity of shocks. As such, the cell averages of the conservative variables throughout the stencil need to be converted to another form prior to reconstruction.

Generically, for a system of $N$ conservation laws, one can fix a transformation matrix $\bm{\Phi}$ of size $N\times N$. Denoting the vector of cell-averaged conservative variables over cell $\Omega_h$ by  $\langle\bm{U}\rangle_{h}$, a conversion from $\langle\bm{U}\rangle_{h}$ to a choice of reconstruction variable $\langle\bm{W}\rangle_{h}$ is done by
\begin{equation}
  \langle\bm{W}\rangle_{h} = \bm{\Phi}\langle\bm{U}\rangle_{h},\quad\forall \Omega_h\in\mathcal{S}_0.
  \label{eq:consToRecVar}
\end{equation}
Reconstruction then proceeds over $\langle\bm{W}\rangle$ component-wise yielding a pointwise estimate of the state, $\bm{w}_{s}$, in the transformed variables. The pointwise defined Riemann state $\bm{u}_{s}$ can be obtained as
\begin{equation}
  \bm{u}_{s}=\bm{\Phi}^{-1}\bm{w}_{s}.
  \label{eq:recVarToCons}
\end{equation}
Crucially, $\bm{\Phi}$ must be constant throughout the stencil for the accuracy of the reconstruction to be maintained. In the following two subsections, we describe two ways to choose the transformation matrix $\bm{\Phi}$. Section \ref{sec:weno:vars:char} briefly overviews the most popular choice with characteristic variables, which is computationally expensive and restricted to governing systems where these variables are known. In Section \ref{sec:weno:vars:prim}, we address these issues by introducing a new approach using the so-called \textit{linearized primitive} variables.

\subsubsection{Characteristic variables}
\label{sec:weno:vars:char}
It is widely known in the computational fluid dynamics community that the best results come from limited reconstruction in characteristic variables, e.g., see \cite{van2006upwind}. For each fixed direction $\bm{\eta}$ and reference state $\langle\widetilde{\bm{U}}\rangle$, the flux Jacobian
\begin{equation}
  \bm{A} = \left.\frac{\partial\bm{F}_{\bm{\eta}}}{\partial\bm{U}}\right|_{\langle\widetilde{\bm{U}}\rangle}
  \label{eq:fluxJac}
\end{equation}
is evaluated from the $\bm{\eta}$ directional flux $\bm{F}_{\bm{\eta}}$. Then the eigendecomposition $\bm{A} = \bm{R}\bm{\Lambda}\bm{L}$ is computed, and the transformation matrices are set as $\bm{\Phi} = \bm{L}$ and $\bm{\Phi}^{-1} = \bm{R}$. The reference state $\langle\widetilde{\bm{U}}\rangle$ is chosen as the cell average in the central cell of the stencil. Note that since the eigenvectors of the flux Jacobian depend nonlinearly on the components of the reference state $\langle\widetilde{\bm{U}}\rangle$, this only provides a second-order accurate approximation of the wave structure within the central cell; however, the limited accuracy of the wave structure does not influence the accuracy of the reconstruction as all that matters is that $\bm{\Phi}$ be fixed over the whole stencil and that $\bm{\Phi}^{-1}$ be its exact inverse. This characteristic reconstruction serves to optimally separate the originally coupled relationships of local flow variables into nearly independent, decoupled characteristic components so that the influence of a discontinuity in one component can be isolated from the other components during reconstruction.

\subsubsection{Linearized primitive variables}
\label{sec:weno:vars:prim}
A critical drawback to using characteristic variables is that the nonlinear weights within the WENO reconstruction must be re-calculated when the direction $\bm{\eta}$ is changed. On structured three-dimensional grids, these nonlinear weights need to be calculated multiple times for each cell, once in the $x-$direction for the reconstruction of the $(i\pm 1/2, j, k)$ Riemann states, again in the $y-$direction for the $(i, j\pm 1/2, k)$ Riemann states, and once again in the $z-$direction for the $(i, j, k\pm 1/2)$ Riemann states. This is of course not an issue for one dimensional reconstruction schemes since they only act in a single fixed direction, but for multidimensional reconstruction it is beneficial to avoid these repeated calculations. The situation worsens further for unstructured grids, with the nonlinear weights being re-evaluated for each and every face of the cell.

A cost-effective alternative can be made available with primitive variables that provide much better reconstruction than conservative variables in the vicinity of shocks with reduced oscillations. By being directionally-independent, they are more computationally friendly than characteristic variables needing only a half or a third as many nonlinear weight calculations in two and three dimensions, respectively (with even greater savings on unstructured grids). Obviously, the primitive variables depend nonlinearly on the conservative variables; converting all cell-averaged conservative variables, $\langle\bm{U}\rangle_h$, in the stencil to primitive variables naively will reduce the accuracy to second-order immediately and irreversibly. To address this, we introduce a new approach below.

Let $\bm{V}$ denote the primitive variables associated with $\bm{U}$, which for the compressible Euler equations are $\bm{V}=\left(\rho,u_1,u_2,u_3,p\right)$, and again set $\langle\widetilde{\bm{U}}\rangle$ as a reference state in conservative variables matching the local cell-average. The irreversible reduction to second-order accuracy can be ameliorated by linearizing the map from conservative to primitive variables around the reference state as
\begin{equation}
  \bm{V}(\bm{U}) \approx \bm{V}(\langle\widetilde{\bm{U}}\rangle) + \left.\frac{\partial\bm{V}}{\partial\bm{U}}\right|_{\langle\widetilde{\bm{U}}\rangle}\left(\bm{U}-\langle\widetilde{\bm{U}}\rangle\right).
  \label{eq:linMap}
\end{equation}
This gives an approximate primitive state $\bm{V}$ corresponding to any given conservative state $\bm{U}$. By linearity, \cref{eq:linMap} can be directly averaged over a cell to obtain an approximation of the average primitive state simply by using a cell-averaged conservative state as the input. This can be applied over the entire stencil to generate approximate cell average values of the primitive variables $\langle\bm{V}\rangle_h$ from each $\langle\bm{U}\rangle_h$. We emphasize that the produced cell averages are only approximate, and do not constitute a high-order representation of the averaged primitive variables.

The constant part of \cref{eq:linMap},
\begin{equation}
\bm{V}(\langle\widetilde{\bm{U}}\rangle) - 
\left.\frac{\partial\bm{V}}{\partial\bm{U}}\right|_{\langle\widetilde{\bm{U}}\rangle}
\langle\widetilde{\bm{U}}\rangle,
\end{equation}
has no influence on the nonlinear weights since it will disappear due to the derivatives in \cref{eq:smoothIndic}. 
Dropping the constant part and applying \cref{eq:consToRecVar}, we get
\begin{equation}
\langle\bm{W}\rangle_h = \left.\frac{\partial\bm{V}}{\partial\bm{U}}\right|_{\langle\widetilde{\bm{U}}\rangle} \langle\bm{U}\rangle_h, 
\quad\forall \Omega_h\in\mathcal{S}_0,
\end{equation}
as the cell-averged values of the linearized primitive variables for reconstruction. Reconstruction can now proceed componentwise over these variables resulting in the pointwise variable $\bm{w}_s$.
Following \cref{eq:consToRecVar} and \cref{eq:recVarToCons}, this establishes the transformation matrices
\begin{align}
 \bm{\Phi} &= \left.\frac{\partial\bm{V}}{\partial\bm{U}}\right|_{\langle\widetilde{\bm{U}}\rangle}, \\
 \bm{\Phi}^{-1} &= \left(\left.\frac{\partial\bm{V}}{\partial\bm{U}}\right|_{\langle\widetilde{\bm{U}}\rangle}\right)^{-1} = \left.\frac{\partial\bm{U}}{\partial\bm{V}}\right|_{\bm{V}(\langle\widetilde{\bm{U}}\rangle)},
 \label{eq:primMat}
\end{align}
where the full forms are given in \cref{eq:lpToConsMat} and \cref{eq:consToLPMat} respectively.
We also note that finding $\bm{\Phi}$ and $\bm{\Phi}^{-1}$ is substantially simpler than finding the eigendecomposition of the flux Jacobian. The entries of both transformation matrices are all available by simply differentiating the maps between the conservative and primitive variables. In particular, the inverse transformation matrix is directly available without needing to perform any symbolic or numerical inversion.

This approach may be interesting for more complicated systems of conservation laws, e.g., the ideal magnetohydrodynamics equations, radiation (magneto)hydrodynamics, multi-physics multi-fluid equations, or systems with non-convex fluxes. As a final note, the transformation matrices themselves contain a large number of zeros and the actual conversions can be done more cheaply by leveraging this fact. This is in contrast to the characteristic variables where the transformation matrices are dense and full.

\subsection{Local smoothness indicators} \label{sec:weno:kxrcf}
The described WENO procedure is effective but considerably more expensive than directly forming an unlimited high accuracy reconstruction using the whole stencil. On the other hand, there are relatively few cells where the nonlinear limited WENO reconstruction is truly needed. For example, the study in \cite{bourgeois2022gp} suggests that, even in highly compressible, shock dominant 2D tests where the initial Mach number reaches as high as 800 (e.g., see Section \ref{sec:results:aj}), only a small fraction of less than 10\% of the entire cells would need nonlinear limiting, while the rest 90\% or more evolve well without any expensive nonlinear limiting.
Hence, if a relatively cheap indicator can flag the cells that actually need WENO, then the cheaper unlimited reconstruction can be used elsewhere. MOOD style schemes (e.g., \cite{clain2011high,bourgeois2022gp}) take this idea to the extreme where unlimited reconstruction is attempted everywhere and order reduction is applied systematically to counteract spurious oscillations near discontinuities. MOOD schemes are inherently \textit{a-posteriori} in their action, which could potentially cause challenging issues when incorporating them into highly parallel solvers for load balancing.

The KXRCF indicator introduced in \cite{krivodonova:KXRCF}, and the related limiter of \cite{fuShu:KXRCF}, operates in an \textit{a-priori} manner as a flag for troubled cells in discontinuous Galerkin (DG) methods. The underlying mechanism for this indicator is related to superconvergence properties of DG methods in the presence of outflows. Fortunately, the overall idea of the KXRCF indicator is general enough to adapt to the present finite volume method, and in fact could potentially be applied more broadly to other \textit{a-priori} limited finite volume methods. Crucially, this indicator is much cheaper to calculate than WENO is to apply, so the added cost of evaluating the KXRCF indicator over all cells is generally negligible. Operationally, if a significant fraction of cells are marked for WENO limited reconstruction, the indicator could be deactivated and WENO reconstruction performed everywhere.

The essential idea behind the KXRCF indicator is to flag cells that exhibit large jumps in some value across cell faces shared with neighboring cells. Unlimited reconstruction is first applied to all cells to generate all Riemann states throughout the domain. Let $s$ index all quadrature points on all faces of a given cell and superscripts $(-)$ and $(+)$ indicate states obtained from the current cell and the relevant adjacent neighbor, respectively. Then the absolute jump of an indicator variable representing some quantity $q$ at the $s^{\rm th}$ quadrature point is $\left|q_s^{(+)} - q_s^{(-)}\right|$. Similarly, let $\langle Q\rangle$ denote the same quantity evaluated using the cell-average. For our purpose here, we use the entropy $q=\frac{p}{\rho^\gamma}$ as the indicator variable, which then yields
$\langle Q\rangle \approx \frac{\langle p \rangle}{\langle \rho \rangle^\gamma}$.

For smooth data the jumps across faces should be small, e.g., $\mathcal{O}\left(\Delta^{2R}\right)$ ideally, while for rough data the jumps should be large. Hence, when the condition for some power $m$
\begin{equation}
  \frac{\max\limits_s \left\{\left|q_s^{(+)} - q_s^{(-)}\right|\right\}}{\langle Q\rangle} < \Delta^m 
  \label{eq:kxrcf}
\end{equation}
holds, the cell does not likely need any form of nonlinear limiting. Alternatively, when the relative jump in $q$ is large at any quadrature point the cell will be flagged as needing WENO. This form is more conservative than that designed in \cite{krivodonova:KXRCF}, but works quite well in the context of the present finite volume method. The power $m$ on $\Delta$ in \cref{eq:kxrcf} controls the sensitivity of the indicator, i.e., larger powers will flag more cells as needing WENO while smaller values will flag less cells permitting larger jumps. For the current study, the power of $m=3/2$ was chosen experimentally as it works well across a wide range of flow conditions and does not appear to need any problem-specific tuning.

\section{Positivity preservation} \label{sec:pos_pres}
The WENO method discussed so far is adequate for many shock problems. However, problems exhibiting very strong shocks or near vacuum states may encounter negative states in the density and/or pressure variables during evolution. The high Mach number astrophysical jets shown in Section \ref{sec:results:aj} are one such example.

To ameliorate issues with negativity, the above method can be combined with the positivity-preserving limiter introduced in \cite{balsara:Astrojets}, which is similar to the methods described in \cite{huAdamsShu:posFVM,zhangShu:posDG}. The innovation in \cite{balsara:Astrojets} is the selection of density and pressure bounds in a data-dependent manner that reduces the impact of extra parameters needing problem-specific tuning.

The current method only differs from \cite{balsara:Astrojets} by doing no reconstructions in the cell interior to test for positivity unless source terms are included in the problem. This reduces the computational cost of the limiter as those interior states would otherwise serve no purpose, but does mean that the limiter is no longer provably positivity-preserving. Within the implementation these additional states can easily be limited by including a zero source term if desired.

A brief description of the limiter as applied in three dimensions will be given here for completeness, and we refer the reader to \cite{balsara:Astrojets} for a more detailed discussion. In what follows let the index $s$ label the individual quadrature points on all faces of the cell in consideration such that $\bm{u}_s$ is the pointwise reconstructed Riemann state at the $s^{\rm th}$ quadrature point. As in Section \ref{sec:weno:vars} let $\langle\widetilde{\bm{U}}\rangle$ be the averaged state in the cell where the limiter is being applied, and let the grid indices $\left(i,j,k\right)$ be defined relative to this cell.

\subsection{Density and pressure bounds} \label{sec:pos_pres:bounds}
The limiter does not enforce any \textit{a-priori} chosen bounds on the density and pressure, e.g., there is no fixed hard floor and/or ceiling values. Instead, these bounds are chosen in a way dependent on the local flow conditions around the cell where limiting is being performed. To do this, we begin with constraining the pointwise density value of each Riemann state from above and below such that $\rho_s\in\left[\rho_{min},\rho_{max}\right]$ with $\rho_{min}>0$, while the pointwise pressure value is constrained only from below such that $p_s\geq p_{min}>0$. 
Using the relative cell index convention for simplicity, these bounds depend upon the local cell average densities and pressures
\begin{align}
  \overline{\rho}_{max} = \max\limits_{-1\leq i,j,k\leq 1}\left\{\langle \rho\rangle_{i,j,k}\right\}, \\
  \overline{\rho}_{min} = \min\limits_{-1\leq i,j,k\leq 1}\left\{\langle\rho\rangle_{i,j,k}\right\}, \\
  \overline{p}_{min} = \min\limits_{-1\leq i,j,k\leq 1}\left\{\langle p \rangle_{i,j,k}\right\},
\end{align}
which are subsequently expanded as 
\begin{align}
  \rho_{max} = \overline{\rho}_{max}\left(1 + \kappa - \kappa\eta\right), \\
  \rho_{min} = \overline{\rho}_{min}\left(1 - \kappa + \kappa\eta\right), \\
   p_{min} = \overline{p}_{min}\left(1 - \kappa + \kappa\eta\right),
\end{align}
where the constant $\kappa$ and the flattener $\eta$ adjust how strict the bounds are by tending to unity only in compressive regions, all of which are set as in \cite{balsara:Astrojets}. Note that the pressures used throughout this section are retrieved directly using cell average values rather than pointwise values and are hence only second-order accurate.

\subsection{Application of the limiter} \label{sec:pos_pres:apply}
With the allowable bounds on density and pressure known, the pointwise Riemann states on the boundary of each cell are modified by the positivity-preserving limiter in two stages as below. The limiter first ensures that the density is bounded above and below such that $\rho_s\in\left[\rho_{min},\rho_{max}\right]$. The size of the needed correction is given by
\begin{equation}
  \theta_\rho = \min\limits_s\left\{1,\frac{\langle\widetilde{\rho}\rangle - \rho_{min}}{\langle\widetilde{\rho}\rangle - \rho_s},\frac{\rho_{max} - \langle\widetilde{\rho}\rangle}{\rho_s - \langle\widetilde{\rho}\rangle}\right\},
  \label{eq:ppThetaRho}
\end{equation}
and the Riemann states are corrected via
\begin{equation}
  \bm{u}_s \leftarrow \langle \widetilde{\bm{U}}\rangle + \theta_{\rho}\left(\bm{u}_s - \langle\widetilde{\bm{U}}\rangle\right).
  \label{eq:ppCorrectDens}
\end{equation}
There are two important observations to make here. First, on a given cell $(i,j,k)$, \textit{all} components of \textit{all} Riemann states are modified in the same way to ensure consistency. Second, there will be no modification to the states if all densities lie in the desired range to begin with.

All Riemann states will now have valid densities, but the pressures may not yet be bounded below by $p_{min}$. As discussed in \cite{balsara:Astrojets,huAdamsShu:posFVM,zhangShu:posDG}, the convexity of the physically admissible region allows this minimum pressure to be enforced again by rewriting the states as a convex combination of their current values and the cell average value $\langle\widetilde{\bm{U}}\rangle$. That is to say,
\begin{equation}
  \bm{u}_s \leftarrow \langle\widetilde{\bm{U}}\rangle + \theta_{p}\left(\bm{u}_s - \langle\widetilde{\bm{U}}\rangle\right),
  \label{eq:ppCorrectPres}
\end{equation}
though now the selection of $\theta_{p}$ is slightly more elaborate than that of $\theta_{\rho}$. In this case, the required correction for the $s^{\rm th}$ quadrature point, denoted as $\theta_{p;s}$, appears as a root of
\begin{equation}
  \label{eq:thetapRoot}
  p\left(\langle\widetilde{\bm{U}}\rangle + \theta_{p}\left(\bm{u}_s - \langle\widetilde{\bm{U}}\rangle\right)\right) - p_{min} = 0.
\end{equation}
As demonstrated in \cite{balsara:Astrojets} this is really a quadratic equation\footnote{See Eq. (12) in \cite{balsara:Astrojets} in which $\tau_{i,j}^q$ corresponds with $\theta_{p;s}$ herein.} in terms of $\theta_{p;s}$ which could be solved analytically. However, for the sake of a uniform implementation, we elect to solve \cref{eq:thetapRoot} using a bisection approach as in the magnetohydrodynamics portion of \cite{balsara:Astrojets}. This allows future inclusions of other equations of state in which \cref{eq:thetapRoot} would fail to be quadratic.

After solving \cref{eq:thetapRoot} for $\theta_{p;s}$ at each quadrature point where the pressure minimum was violated, the overall correction is selected as the smallest admissible root
\begin{equation}
  \label{eq:thetaP}
  \theta_p = \min\limits_s\left\{\theta_{p;s} ~|~ 0\leq\theta_{p;s}\leq 1\right\},
\end{equation}
and \cref{eq:ppCorrectPres} is applied. As before, it should be noted that all components of all Riemann states are modified in unison, and no correction will be applied if all pressures appearing in the Riemann states already satisfy $p_s\geq p_{min}$. Additionally, this correction naturally leaves the densities bounded appropriately. As a final note, the presentation here only considers the Riemann states on the boundaries of a cell, but interior states could easily be included in the set $\lbrace \bm{u}_s\rbrace$.

\section{Overview and implementation} \label{sec:method}
All of the building blocks of the proposed scheme have been defined, and now we gather them together into an explicit step-by-step description. The application of the solver proceeds in three main stages, where the first and the second stages are conducted one-time as part of the initial setup, after which the last stage is performed repeatedly to evolve solutions until the final simulation time. The three main stages include:
\begin{enumerate}
\item[(i)] constructing the stencils and reconstruction vectors, 
\item[(ii)] evaluating initial conditions into cell average quantities, and finally, 
\item[(iii)] advancing the cell averages through time. 
\end{enumerate}
We elaborate on each of these stages in the following subsections.

\subsection{Stencils and weight vectors} \label{sec:method:precomp}
To prepare the stencils for use in each uniform grid simulation, one needs to:
\begin{enumerate}
\item Enumerate all cells in the full stencil, $\mathcal{S}_0$, and the remaining substencils, $\mathcal{S}_q,~q=1,\ldots,N_S$.
\item Form reconstruction vectors, $\bm{r}_{(q)}$, for each face quadrature point $\bm{x}_s$ relative to each (sub)stencil using \cref{eq:recVec1,eq:recVec2,eq:recVec3}. The quadrature points on each face are set using a Gauss-Legendre quadrature rule with $R$ points per dimension. This yields $2R^{\rm th}$ order accuracy matching that of the reconstruction scheme.
\item Form reconstruction vectors for all partial derivatives needed to evaluate the smoothness indicators in \cref{eq:smoothIndic} relative to each (sub)stencil.
\end{enumerate}
The choice of the hyperparameter $\ell$ in the kernel \cref{eq:kernSE} deserves further discussion. Kernel-based interpolation based on the squared exponential kernel can be used without the trailing monomial terms in \cref{eq:avExpan}. In this case maintaining a consistent order of accuracy requires that the hyperparameter $\ell$ be fixed independent of the grid scale. As a result, the kernel matrix in \cref{eq:recBlockSys} becomes increasingly ill-conditioned as the grid is refined, and a number of elaborate stable inversion schemes have been developed to handle this issue \citep{fornberg:rbfpade,wright2017stable,fornberg:rbfga}. 

Instead, another resolution is made available by including the monomial terms in \cref{eq:avExpan}, which yields guaranteed convergence rates independent of $\ell$.  Larger values of $\ell$ are desirable as discussed in \cite{flyer:RbfWithPolys}. Then, a good way to proceed is to set the monomial degree as high as possible, whereby one can increase $\ell$ as far as allowed without needing to resort to the more elaborate stable inversion schemes that avoid (near) singular matrix inversions \citep{fornberg:rbfpade,wright2017stable,fornberg:rbfga}. This approach with the monomial terms is the strategy we take here, and over all of the various stencils we have found that $\ell = 5\Delta$ is a good baseline choice. Larger values are possible, but none of the presented results are particularly sensitive to this choice.

\subsection{Evaluation of initial conditions} \label{sec:method:ics}
The initial conditions must be set accurately if the succeeding evolution is to be meaningful, e.g., see \cite{bourgeois2022gp}. Here the initial cell averages are filled using a tensor-product Gauss-Legendre rule on each cell. As for the quadrature points on faces, we use a rule with $R$ points per dimension (see bullet point 2 in Section \ref{sec:method:precomp}).

\subsection{Time advancement} \label{sec:method:time}
We describe the overall solution advancement of the KFVM-WENO solver whose high-order spatial solutions are temporally evolved by multistage Runge-Kutta (RK) type integrators. The implementation of these time integrators ultimately comes down to being able to evaluate all spatial terms with respect to some given state. The evaluation of the spatial terms proceeds as follows:
\begin{enumerate}
\item Fill the ghost cells in accordance with the boundary conditions.
\item Calculate Riemann states at each quadrature point on cell faces (see below).
\item Apply the positivity preserving limiter on the Riemann states, following Section \ref{sec:pos_pres}.
\item Populate Riemann states outside of the physical domain in accordance with the boundary conditions.
\item Call an approximate Riemann solver at each flux quadrature point.
\item Integrate fluxes and update cell averages in accordance with the chosen RK method.
\end{enumerate}
Naturally, the second step of calculating the high-order KFVM-WENO Riemann states contains most of the contributions from this work. This proceeds as follows:
\begin{enumerate}
\item Reconstruct unlimited Riemann states on all cells using cell averages from the full stencils $\mathcal{S}_0$.
\item Evaluate the KXRCF indicator from Section \ref{sec:weno:kxrcf}, and flag cells needing WENO.
\item For all flagged cells do the following:
  \begin{enumerate}
  \item Construct transformation matrix $\bm{\Phi}$ from \cref{eq:primMat} using the central cell average data as the reference state.
  \item Project the average state for each cell in $\mathcal{S}_0$ onto the linearized primitive variables via \cref{eq:consToRecVar}.
  \item Reconstruct pointwise values of linearized primitive variables at all face quadrature points via the KFVM-WENO reconstruction detailed in Section \ref{sec:weno:ao}.
  \item Project the pointwise linearized primitive values back to conservative variables via \cref{eq:recVarToCons}, yielding the Riemann states.
  \end{enumerate}
\end{enumerate}

\subsection{Parallel code implementation} \label{sec:code_implementation}
The presented scheme has been implemented using Kokkos \citep{kokkosMain,kokkosOriginal} for shared memory parallelism combined with the message passing interface (MPI) for distributed memory parallelism. This allows for straightforward use of multiple GPUs for high degrees of parallelism. Internally, the reconstruction process makes extensive use of the batched linear algebra routines from the Kokkos-Kernels package \citep{kokkosEcosystem}. We encourage the interested reader to consult our code for more details \citep{may:KFVMKokkos}.

\section{Numerical results} \label{sec:results}
The proposed scheme is evaluated against a variety of benchmark problems, several of which are known to be very stringent in the literature. Below, we briefly overview what is expected in each test problem in this section.

First we experimentally verify the expected orders of accuracy through a convergence study using the isentropic vortex problem \citep{shu1998essentially,spiegel2015survey} in Section \ref{sec:results:isenvort}. The Sod shock tube problem \citep{sod1978survey} is then solved in grid-aligned and tilted configurations in Section \ref{sec:results:sod} as a first evaluation of the shock handling capabilities, and to test if the method has any inherent preference for grid-aligned phenomena.

With those two fundamental tests done we move to more challenging problems. Section \ref{sec:results:rm} presents a Richtmeyer-Meshkov instability driven by a Mach 3 shock wave as described in \cite{samtaney:rmi}. Next, we solve two very stringent problems in Section \ref{sec:results:aj}, where we consider the high Mach astrophysical jet problems from \cite{balsara:Astrojets}.

Finally, two problems with physical viscosity are solved. The Taylor-Green vortex is shown in Section \ref{sec:results:tg} and compared against the validated benchmark data from \cite{Nasa_CFD_workshop_TG}, and Section \ref{sec:results:rt} considers a viscous version of the Rayleigh-Taylor instability as described in \cite{shiRayTay:2003}.

All problems are solved using explicit Runge-Kutta time integrators with error-estimate-based step size selection in addition to the more standard CFL constrained step sizes. To this end we utilize the $[3S]_+^*$ and $[4S]_+^*$ time integrators from \cite{ranocha:optRkErrEst}, which also provides an excellent discussion of error-based step size selection. Due to our treatment of problems with very strong shocks, we also consider the RK-SSP(4,3) method of \cite{kraaijevanger:ssp43} with the embedded RK-SSP(3,2) method described in \cite{fekete:sspEmbMeth}, the combination of which is hereafter denoted as SSP(4,3,2).

The flux integrals are evaluated using tensor-product Gauss-Legendre quadrature rules with $R$ points per dimension as discussed in Section \ref{sec:weno}. The length scale hyperparameter in \cref{eq:kernSE} is always set as $\ell = 5\Delta$. All problems use the newly proposed linearized primitive variables from Section \ref{sec:weno:vars:prim}. As for Riemann problems, the HLLC+ approximate Riemann solver from \cite{chen:hllcPlus} is applied to all problems apart from the astrophysical jets which use the more stable HLL \citep{hartenLaxLeer:HLL,toro} approximate Riemann solver. In each case wavespeeds are estimated following \cite{batten:wavespeed}.

\begin{table}[h]
  \caption{Shown are the experimental orders of convergence (EOC) for the described method 
  with the radius $R=2$ scheme (left) and the $R=3$ scheme (right) as tested on the isentropic vortex problem.}
  \label{tab:ivConvR2nwR3nw}
{\centering
  \begin{tabular}{@{}r|cc|cccc@{}}
    \toprule
\multicolumn{1}{c}{Grid} & \multicolumn{2}{|c}{$R=2$} & \multicolumn{2}{|c}{$R=3$} \\
\cmidrule(rl){2-3} \cmidrule(rl){4-5}
    Res.        & $L_1$ Error    & EOC                        & $L_1$ Error    & EOC       \\ \midrule
    $32^2$    & $1.45\times 10^{-3}$ & --                   & $7.51\times 10^{-4}$ & --                    \\  
    $64^2$    & $2.27\times 10^{-4}$ & \textbf{2.68} & $7.50\times 10^{-5}$ & \textbf{3.32}  \\
    $128^2$   & $1.47\times 10^{-5}$ & \textbf{3.95} & $1.46\times 10^{-6}$ & \textbf{5.68} \\
    $256^2$   & $5.39\times 10^{-7}$ & \textbf{4.77} & $1.76\times 10^{-8}$ & \textbf{6.38} \\ \bottomrule
  \end{tabular}\par}
\end{table}

%

\subsection{Isentropic vortex} \label{sec:results:isenvort}
The isentropic vortex is one of the few nonlinear problems in the literature where its initial condition serves as a smooth, exact solution to the compressible Euler equations. By design, the isentropic vortex evolves in a fully nonlinear manner. These features make it an ideal candidate for experimentally validating the convergence rate of a given method. Herein we set up the problem as discussed in \cite{spiegel2015survey}, which uses a domain of $\Omega = [-10,10]^2$ with periodic boundary conditions and initial conditions given by
\begin{align}
  \rho &= \left(1 + \frac{1-\gamma}{2}\omega^2\right)^{\frac{1}{\gamma - 1}}, \\
  u &= 1 - y\omega, \\
  v &= 1 + x\omega, \\
  w &= 0, \\
  p &= \frac{1}{\gamma}\left(1 + \frac{1-\gamma}{2}\omega^2\right)^{\frac{\gamma}{\gamma - 1}},
\end{align}
where $r=1$ is the vortex radius and the rotation rate $\omega$ is set as
\begin{equation}
  \omega = 5\frac{\sqrt{2e}}{4\pi} e^{-\frac{1}{2}(x^2+y^2)}.
\end{equation}

After evolving to the final time $t=20$, the vortex returns to its initial position, and the accumulated error during the run can be found by comparing the initial and final states. 
In \cref{tab:ivConvR2nwR3nw}, 
the experimental $L_1$ convergence rates of the method are demonstrated for radius $R=2$ and $R=3$ stencils respectively. In all cases the $[4S]^*_+$ time integrator \citep{ranocha:optRkErrEst} is used with the tolerances set as $atol=rtol=10^{-6}$ which dominate over a maximal CFL of 1.0. This tolerance is tight enough to ensure that the spatial error dominates in every case. In \cref{tab:ivConvR2nwR3nw}, we display the obtained $L_1$ errors and associated experimental order of convergence (EOC) rates
\begin{equation}
\text{EOC} = \frac{\ln \left(E_c/E_r\right)}{\ln(2)},
\end{equation}
where $E_c$ and $E_r$ are the errors in the $L_1$ norm on the coarse and the next coarse resolutions (e.g., if $E_c$ is measured on $32^2$, $E_r$ is on $64^2$)
%
for the radius $R=2$ and $R=3$ schemes. 

For radius $R=2$ we expect at least $4^{\rm th}$-order accuracy. Indeed, for resolutions beyond $64\times 64$ we see convergence rates matching or exceeding the $4^{\rm th}$-order convergence. As noted in Section \ref{sec:recon:asym} the polynomial tail in \cref{eq:avExpan} is not necessary to obtain a well-behaved and accurate scheme, but is rather present to allow smaller choices for $\ell$ and a better conditioned linear system in \cref{eq:recVec1,eq:recVec2,eq:recVec3}. Indeed the kernel part of the expansion still plays an active role in improving the accuracy of the overall scheme, and better than $4^{\rm th}$-order accuracy is observed in \cref{tab:ivConvR2nwR3nw}. Similarly, for radius $R=3$, we expect at least $6^{\rm th}$-order accuracy and again we observe orders near or beyond the expected rate after a resolution of $64\times 64$.

\subsection{Sod shock tube} \label{sec:results:sod}
The Sod shock tube problem \citep{sod1978survey} is staple benchmark for any numerical method treating the compressible Euler equations. We consider this problem in two configurations. First, in the standard grid-aligned formulation we take the domain to be $\Omega=[0,1]\times[0,0.04]$ with outflow conditions in the $x-$direction and periodic conditions in the $y-$direction. The initial conditions are piecewise constant with $\rho_L=1$ and $p_L=1$ set for $x<0.5$ and $\rho_R=0.125$ and $p_R=0.1$ otherwise, and zero velocity throughout.

Second, we consider a tilted configuration where the same solution evolves oblique to the grid 
by the angle of $\theta \approx {26.5651^{\circ}}$ (or $\tan \theta = 1/2$) to test if there is any preference for grid-aligned phenomena. Following the idea in \cite{kawai2013divergence,lee2021recursive}, the domain is set as $\Omega=[0,\sqrt{5}]\times[0,2\sqrt{5}]$ with periodic conditions in both directions. The tilted coordinate is set as $x_\parallel = \frac{1}{\sqrt{5}}\left(2x + y\right)$ with $0 \le x_\parallel \le 4$, and the \textit{left} state is initialized in the regions $0 \le x_\parallel \le 0.5$, $1.5 < x_\parallel \le 2.5$, and $3.5 < x_\parallel \le 4$, while the \textit{right} state is imposed elsewhere. In this configuration the final density and pressure extracted along the line $0 \leq x_\parallel \leq 1$ will match that of the grid aligned configuration.

\begin{figure}
  \centering
  \includegraphics[width=\linewidth]{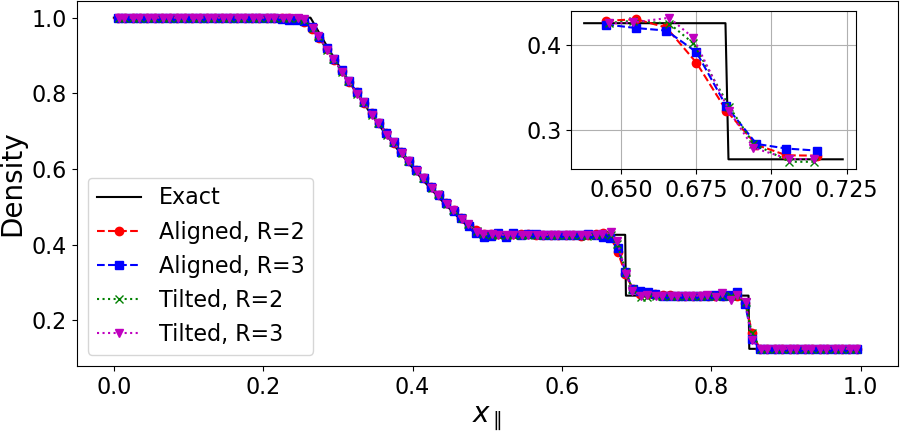}
  \caption{Shown is a trace of the density in the Sod shock tube problem as obtained from four different cases. The solid black line shows the exact analytical solution. The dashed lines with circle and square markers show the results of the radius $R=2$ and $R=3$ schemes in the grid-aligned configuration. The dotted lines with cross and triangle markers show the results of the radius $R=2$ and $R=3$ schemes in the grid-tilted configuration. The inset shows a zoom-in of the region near contact discontinuity.}
  \label{fig:sodtilt}
\end{figure}

These problems are solved with both the radius $R=2$ and $R=3$ schemes to the final time $t=0.2$ using the SSP(4,3,2) time integrator with tolerances $atol=rtol=10^{-4}$ and a maximum CFL of 1.25. For the grid-aligned configurations a grid with spacing $\Delta=1/100$ is used. For the tilted configuration a comparable resolution is achieved by using a grid with spacing $\Delta=\sqrt{5}/250$, which yields 100 cells along the line $0 \leq x_\parallel \leq 1$ where the solution is extracted. The resulting density traces can be seen in \cref{fig:sodtilt}, and all schemes can be seen to produce excellent results in comparison to the exact solution. The inset shows that all schemes resolve the contact discontinuity in roughly three cells with only a small overshoot from the $R=3$ scheme in the tilted configuration. Overall, the present scheme has a minimal preference between features irrespective of the underlying grid configurations.

\subsection{Richtmeyer-Meshkov instability} \label{sec:results:rm}
\begin{figure}
  \centering
  \includegraphics[width=\linewidth,trim={0 1.5in 0 0},clip]{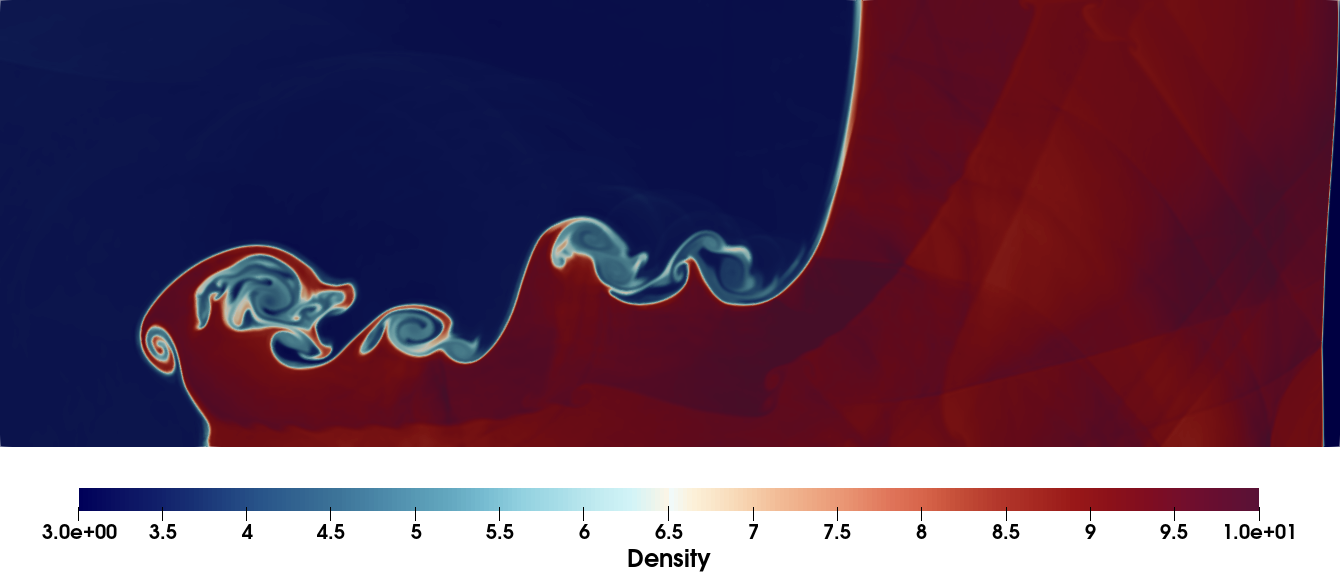}
  \includegraphics[width=\linewidth]{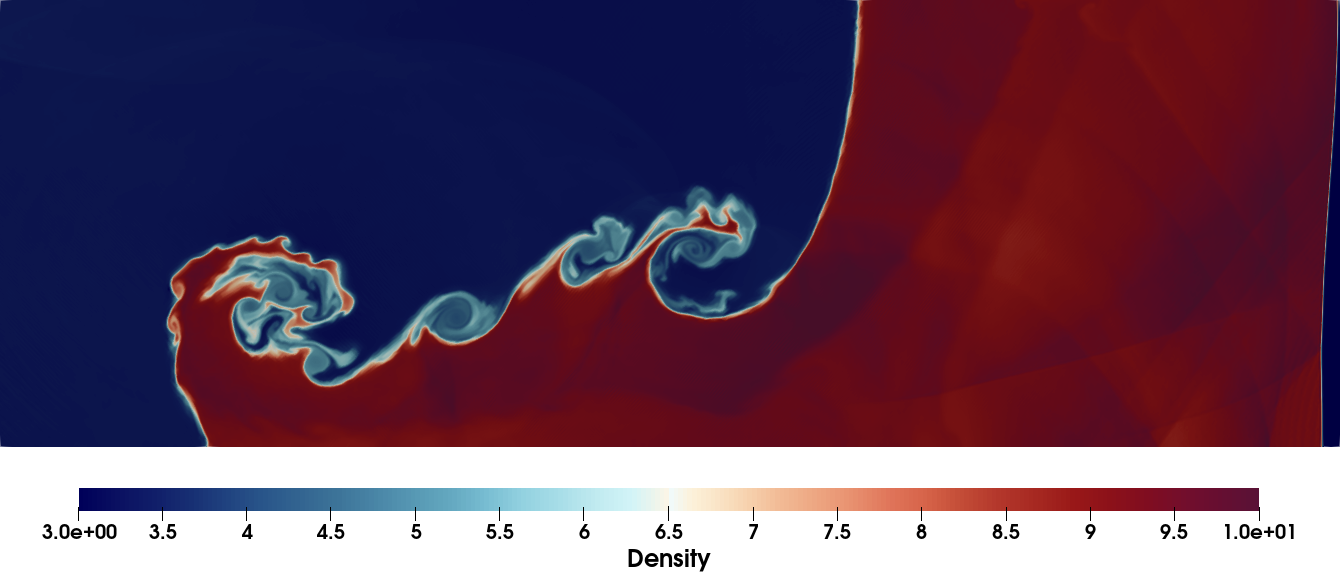}
  \caption{Shown is the density field for the Richtmeyer-Meshkov instability at the final time of $t=3.33$ as solved by the radius $R=2$ scheme (top) and the radius $R=3$ scheme (bottom) on a grid with spacing $\Delta=1/512$. We display the view zoomed into the region $[5/2,11/2]\times [0,1]$ to highlight the interface.}
  \label{fig:rmiR2R3}
\end{figure}

Here we consider a Richtmeyer-Meshkov instability similar to the unmagnetized case in \cite{samtaney:rmi}. A right-traveling shock impinges on an oblique density jump, and the initially straight interface is bent as the shock diffracts through it. As time progresses secondary Kelvin-Helmholtz instabilities are excited along the density interface. The domain is $\Omega = [-1/2,11/2]\times[0,1]$ with fixed inflow conditions at $x=-1/2$, extrapolation conditions at $x=11/2$ and reflecting conditions in the $y-$direction. The problem is parameterized by the shock Mach number $Ma$ and the density $\rho_D$ of the gas downstream of the interface. The initial density is given by
\begin{equation}
  \rho = \begin{cases} \left(1 - \frac{2}{\gamma + 1}\left(1 - \frac{1}{Ma^2}\right)\right)^{-1},~&x<0.2, \\ 1,~&x<y, \\ \rho_D,~&x\geq y, \end{cases}
\end{equation}
the initial pressure is given by
\begin{equation}
  p = \begin{cases} 1 + \frac{2\gamma}{\gamma + 1}\left(Ma^2 - 1\right),~&x<0.2, \\ 1,~&x\geq 0.2, \end{cases}
\end{equation}
and finally the initial $x-$velocity is given by
\begin{equation}
  u = \begin{cases} Ma\sqrt{\gamma}\left(1 - \frac{1}{\rho}\right),~&x<0.2, \\ 0,~&x\geq 0.2. \end{cases}
\end{equation}
As in \cite{samtaney:rmi} we set the parameters $Ma=3$ and $\rho_D=2$. 
In \cref{fig:rmiR2R3}, we show the density fields around the interface at the final time of $t=3.33$ on grids with spacing $\Delta=1/512$ as produced by the radius $R=2$ and $R=3$ schemes, respectively. In both cases the SSP(4,3,2) time integrator is used with tolerance $atol=rtol=10^{-3}$ and a maximum CFL of $1.0$.

The density profiles shown in \cref{fig:rmiR2R3} qualitatively agree with the results presented in \cite{samtaney:rmi} despite the latter making extensive use of adaptive mesh refinement to reach an effective grid resolution of $16384\times 2048$, compared to our results on $1536 \times 512$. The radius $R=3$ results presented here exhibit smaller scale structures along the primary density interface separating the heavy and light fluids than are present in the radius $R=2$ results. The density interfaces within the heavy fluid that arise from the reflections of passing shock are similar in all cases. Both results also pick up an additional Kelvin-Helmholtz instability along an internal interface below the primary one, more closely matching the referenced results despite using a grid with four times larger spacing.

\begin{figure}
  \centering
  \includegraphics[width=\linewidth]{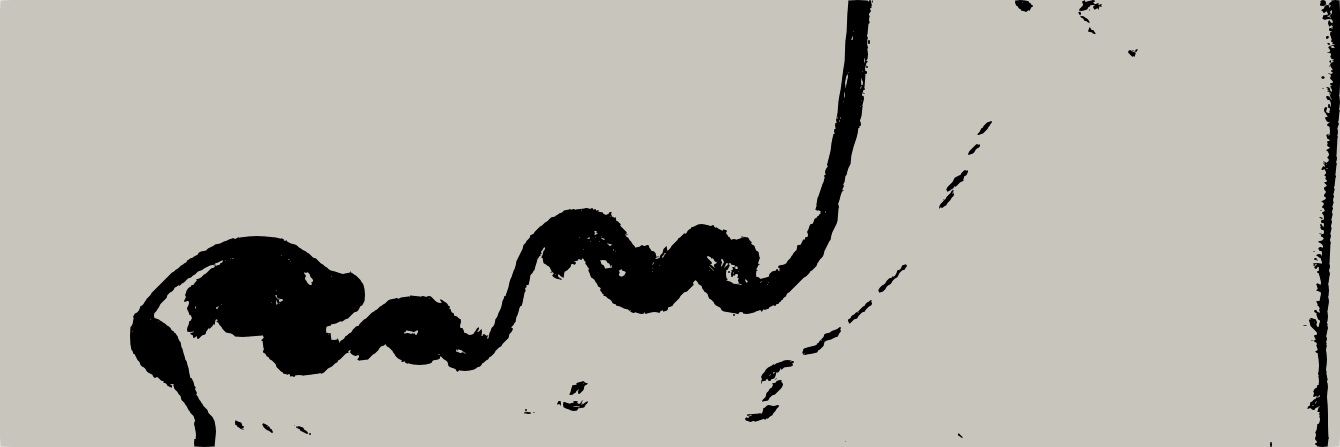}
  \caption{Shown in black are the cells flagged for WENO reconstruction at the final time for the radius $R=2$ scheme. The view has been zoomed in to match \cref{fig:rmiR2R3}.}
  \label{fig:rmiWenoFlag}
\end{figure}

This test also showcases the benefits of the KXRCF style indicator described in Section \ref{sec:weno:kxrcf}. The cells flagged as needing WENO reconstruction in the final time step are shown in \cref{fig:rmiWenoFlag}, and comprise only $5.56\%$ of the overall grid. This problem was run on an NVidia 4080 with the Kokkos profiling tools. With the indicator enabled and WENO applied only sparsely, the solver spent a total of $117$ seconds doing linear reconstruction and $56$ seconds doing WENO reconstruction over the entire evolution of the problem, which consisted of $21{,}193$ right hand side evaluations. With the indicator disabled the solver always performs WENO reconstruction, which took $690$ total seconds across $21{,}349$ right hand side evaluations.

We also note that the cost of activating KXRCF is negligible. The average time to execute one right hand side evaluation (which includes all reconstruction, limitation, Riemann solves, etc.), with the indicator active was $0.044$ seconds and without the indicator the average time was $0.050$ seconds. The overhead in applying the indicator and doing a two-pass reconstruction is easily worthwile, and would only become more useful in three dimensions or with larger stencils.

\subsection{Astrophysical jets} \label{sec:results:aj}
Here we consider the high-density and low-density astrophysical jet problems from \cite{balsara:Astrojets,bourgeois2022gp,ha2008positive} as well as an extension of the low-density jet to three dimensions. The domain used is $\Omega=[0,1/2]\times[0,3/2]$ in two dimensions and $\Omega=[0,1/2]\times[0,3/2]\times[0,1/2]$ in three dimensions. Reflecting conditions are applied on the lower boundaries of $x$ (for 2D and 3D) and $z$ (for 3D only), outflow conditions are applied on the upper boundaries of $x$, $y$ (for 2D and 3D) and $z$ (for 3D only), and mixed inflow/outflow conditions are applied on the lower $y$ boundary. The highly compressible extragalactic jet enters the ambient flow through a narrow slit configured as $r^2 = x^2 < 0.05^2$ for 2D and $r^2 = x^2+z^2 < 0.05^2$ for 3D, where the inflow conditions are imposed on the slit. Apart from the narrow slit, the remainder of the lower $y$ boundary is set to simple outflow extrapolation conditions. 

In all cases, the injected jet is initialized with density $\rho=\gamma$ and unit pressure $p=1$ and all runs are solved using the HLL approximate Riemann solver and the SSP(4,3,2) time integrator with tolerances of $atol=rtol=10^{-2}$ and a maximum CFL of $1.0$.

In the high-density case, the jet evolves with a $y-$velocity of $800$ into an initially quiescent background with density $\rho=\gamma/10$ and unit pressure, making the jet's density is ten times higher than the background density. On the other hand, in the low-density case, the jet is ten times lighter than the background ambient flow that has density $\rho=10\gamma$ and unit pressure, moving with a $y-$velocity of $100$. All other velocity components are zero in all cases. These configurations yield jet Mach numbers of $800$ and $100$, respectively for the high-density jet and the low-density jet cases.

\begin{figure}
  \centering
  \includegraphics[width=\linewidth]{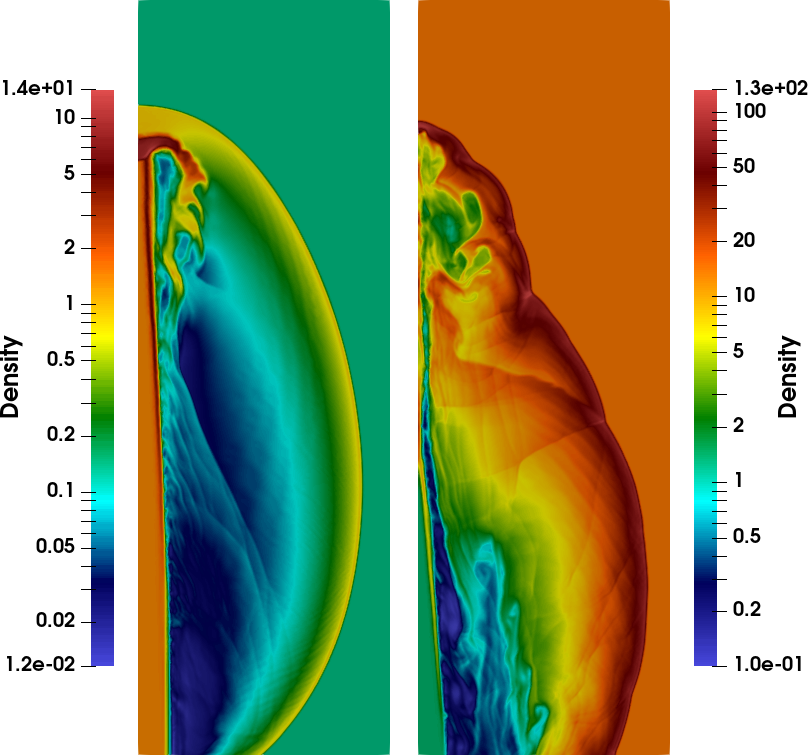}
  \caption{Shown are the logarithmic density fields for the high-density (left) and low density (right) astrophysical jets in two dimensions at the final times of $t=0.002$ and $t=0.04$ respectively. In both cases radius $R=3$ stencils were used on grids with spacing $\Delta=1/512$.}
  \label{fig:astrojet2D}
\end{figure}

\begin{figure}
  \centering
  \includegraphics[width=\linewidth]{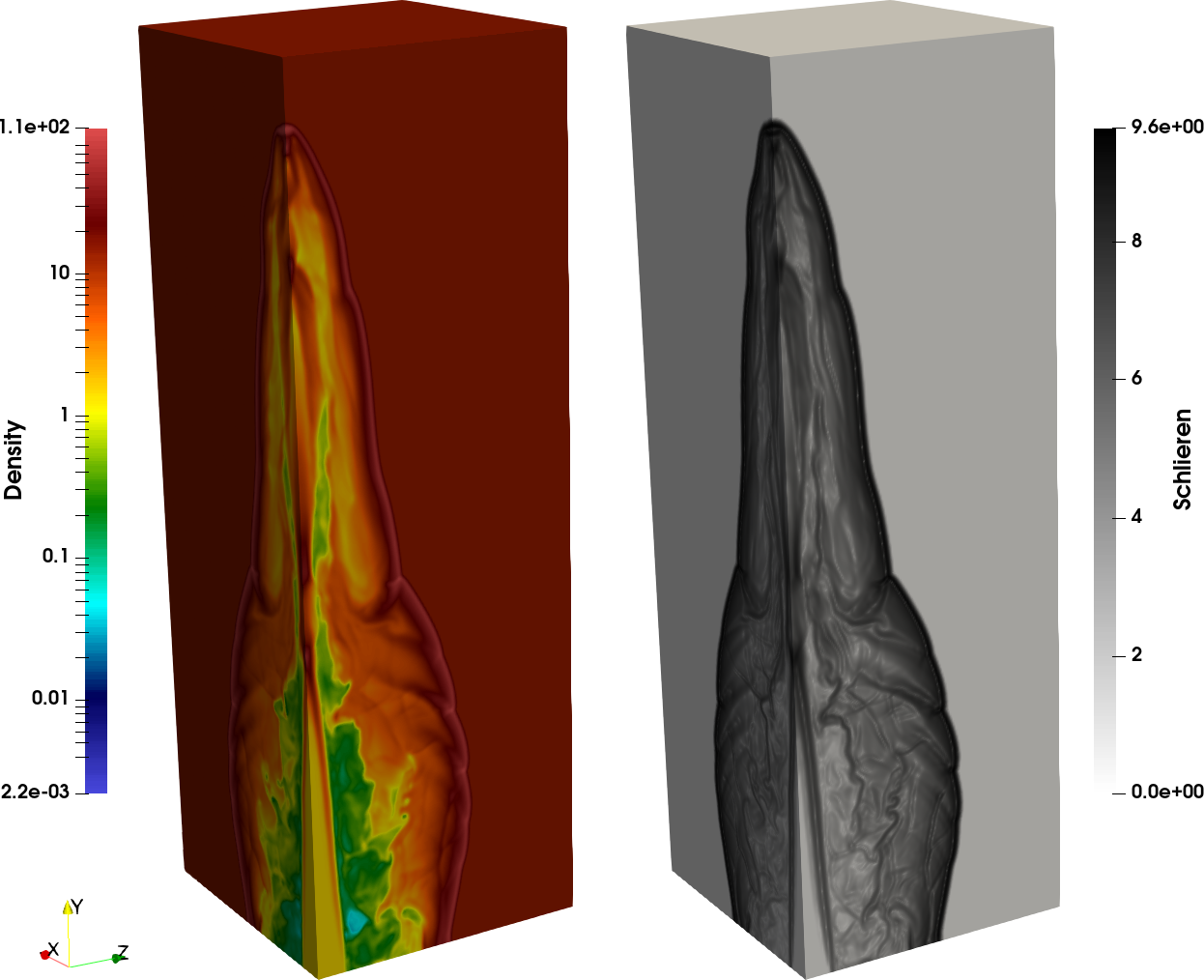}
  \caption{The left panel shows the logarithmic density field for the low-density astrophysical jet in three dimensions at the final time of $t=0.035$ as solved by the radius $R=2$ scheme on a grid with spacing $\Delta=1/384$. The right panel shows the corresponding numerical Schlieren image defined as $\ln\left(1+\left|\nabla\rho\right|\right)$.}
  \label{fig:astrojet3D}
\end{figure}

In two dimensions both the high and low density jets are solved using the radius $R=3$ scheme to final times of $t=0.002$ and $t=0.04$ respectively on grids with spacing $\Delta=1/512$, and the resulting logarithmic density fields can be seen in \cref{fig:astrojet2D}. In three dimensions we consider only the more stringent low-density jet using the radius $R=2$ scheme on a grid with spacing $\Delta=1/384$. The logarithmic density field at the final time $t=0.035$ can be seen in \cref{fig:astrojet3D}.

As seen in \cref{fig:astrojet2D} and  \cref{fig:astrojet3D}, the overall evolution of the jets, the flow dynamics surrounding the jet envelops, the internal flow structures such as Kelvin-Helmholtz instabilities in the ``cocoon'' region, are all well captured and show excellent agreement with the results studied in \cite{balsara:Astrojets,bourgeois2022gp,ha2008positive}. Our results assure that the present schemes are capable of simulating highly compressible flows where the shock strengths are much beyond the typical supersonic, hypersonic, and high-hypersonic regimes.

\subsection{Taylor-Green vortex} \label{sec:results:tg}
The Taylor-Green vortex studied in \cite{taylor1937mechanism} furnishes a prototypical transition to turbulence and turbulence decay problem. Initially, smooth large-scale vortices decay into smaller and smaller vortices and eventually fully decay. The Taylor-Green vortex is well studied in literature especially with respect to high order methods and is therefore an ideal test case for our schemes. 

We consider the Taylor-Green vortex at Reynolds number $Re=1{,}600$ following its prescription as a challenge problem reported using incompressible flow solvers (e.g., see \cite{van2011comparison} for a comparison study using a particle-mesh vortex method and a pseudo-spectral method) as well as in several code-to-code comparison workshops on high-order methods for computational fluid dynamics. 

Reference scales $\widetilde{\rho},U,L$ are fixed for the density, velocity, and length respectively. The induced time scale is $\tau=L/U$, and the pressure scale is chosen to be $P=\widetilde{\rho}U^2$. After nondimensionalizing relative to these reference scales, we solve the problem in a triply periodic box on the domain $\Omega=[-\pi,\pi]^3$, with the initial conditions
\begin{align}
  \rho &= 1 + \frac{\gamma Ma^2}{16}\left(\cos(2x)+\cos(2y)\right)\left(\cos(2z)+2\right), \nonumber \\
  u &= \sin(x)\cos(y)\cos(z), \nonumber \\
  v &= -\cos(x)\sin(y)\cos(z), \nonumber \\
  w &= 0, \nonumber \\
  p &= \frac{1}{\gamma Ma^2} + \frac{1}{16}\left(\cos(2x)+\cos(2y)\right)\left(\cos(2z)+2\right), \nonumber
\end{align}
which correspond to constant initial temperature as $\rho = \gamma Ma^2 p$. We follow the guidance on the problem statement found in one of the CFD workshops held at the NASA Glenn Research Center
\footnote{See for instance the problem description of case C3.3 at \url{https://www1.grc.nasa.gov/research-and-engineering/hiocfd/}} to set
the Mach number as low as $Ma=0.1$ to approximate incompressible flow. The nondimensional governing equations are
\begin{equation}
  \frac{\partial \bm{U}}{\partial t} + \nabla\cdot\left(\bm{F} - \bm{F}^{(v)}\right) = 0,
  \label{eq:ns}
\end{equation}
where the vector of conserved quantities $\bm{U}$ and the inviscid fluxes $\bm{F}$ are those in \cref{eq:euler}. The additional viscous fluxes are given by
\begin{equation}
  \bm{F}^{(v)}_j = \begin{pmatrix} 0 \\ \sigma_{ij} \\ q_j + \sigma_{jk}u_k \end{pmatrix},
  \label{eq:nsVisc}
\end{equation}
where the scaled stress tensor is
\begin{equation}
  \sigma_{ij} = \frac{1}{Re}\left(\frac{\partial u_i}{\partial x_j} + \frac{\partial u_j}{\partial x_i} - \frac{2}{3}(\nabla\cdot\bm{u})\delta_{ij}\right),
  \label{eq:shearTensor}
\end{equation}
and the thermal conduction is
\begin{equation}
  q_j = \frac{1}{PrRe}\frac{\partial T}{\partial x_j},
  \label{eq:thermalFlux}
\end{equation}
where the Prandtl number has been fixed at $Pr=0.71$. Since this paper has the compressible Euler equations as the primary focus we elect to discretize these viscous terms through the standard second order accurate finite differences for simplicity to handle viscous effects. This is justified in this high Reynolds number case, although future extensions to high accuracy viscous terms would be interesting.

\begin{figure}
  \centering
  \includegraphics[width=\linewidth]{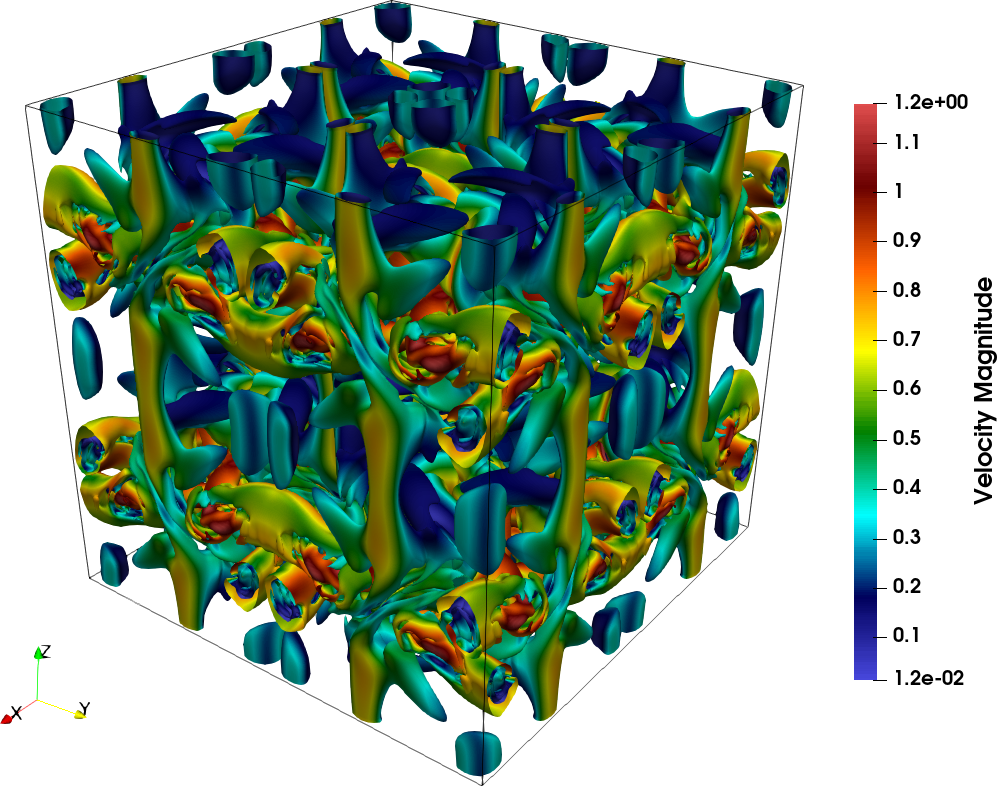}
  \caption{The $Q=0.1$ isosurface of the Q-criterion colored by velocity magnitude is shown for the Taylor-Green vortex at non-dimensional time $t=8$ on a grid with spacing $\Delta=\pi/192$.}
  \label{fig:tgVort}
\end{figure}

This problem is solved using a radius $R=2$ stencil on a grid with spacing $\Delta = \pi/192$ yielding $384^3$ cells in total. The $Q=0.1$ isosurface of the $Q-$criterion colored by velocity magnitude at non-dimensional time $t=8$ is shown in \cref{fig:tgVort}. Our result can be compared directly with the results presented in \cite{giangaspero2015}. 

Running the Taylor-Green vortex problem is considered a challenging benchmark code-to-code verification test for a wide range of high-order schemes. As is standard for this challenge problem, the fully converged high resolution data available from \cite{Nasa_CFD_workshop_TG} is used to provide quantitative validation of the present scheme. These results were obtained using a pseudo-spectral scheme for the incompressible Navier-Stokes equations on a grid of size $512^3$, the details of which can be found in the associated problem statement.

\begin{figure}
  \centering
  \includegraphics[width=\linewidth]{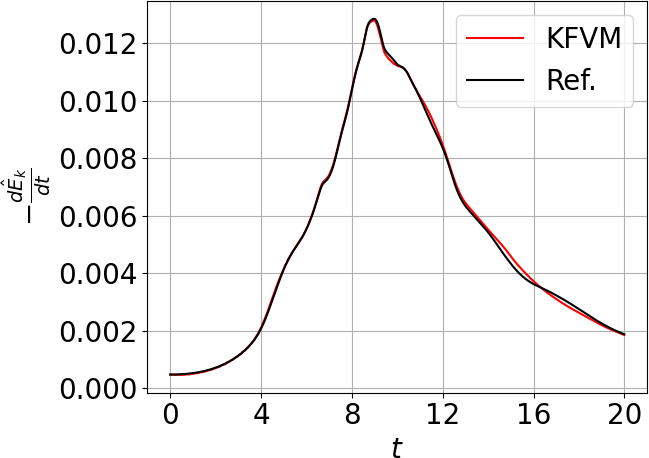}
  \caption{The evolution of the kinetic energy dissipation rate is compared against the fully converged reference data from \cite{Nasa_CFD_workshop_TG}.}
  \label{fig:tgVortStats}
\end{figure}

In \cref{fig:tgVortStats} we compare the (non-dimensional) kinetic energy dissipation rate, $-\frac{d\hat{E_k}}{dt}$, from the present scheme to the reference data. We observe excellent agreement between our scheme and the reference data through the peak dissipation rate at time $t\approx 9$, with only minor discrepancies through the decay process after the peak.

\subsection{Rayleigh-Taylor instability} \label{sec:results:rt}
Here we consider a two dimensional viscous Rayleigh-Taylor instability. The domain and initial conditions follow the setup from \cite{shiRayTay:2003}, though now an additional physical viscosity is added. The governing equations follow from the previous section (see \cref{eq:ns}) with an added source term in \cref{eq:nsVisc}.

Gravity is taken to point in the positive $y-$direction, and nondimensionalizing as before yields the source term
\begin{equation}
  \label{eq:rtSource}
  \bm{S} = \frac{1}{Fr^2}\begin{pmatrix} 0, & 0, & \rho, & 0, & \rho u_2 \end{pmatrix}^T,
\end{equation}
where $Fr=\frac{U}{\sqrt{gL}}$ is the Froude number. We set the domain as $\Omega = [0,1/4]\times[0,1]$. Initially, the density is set to $\rho_{high}=2$ for $y<1/2$ and $\rho_{low}=1$ otherwise, the pressure is set as
\begin{equation*}
  p = \begin{cases} \frac{\rho_{high}}{Fr^2}y+ 1,\quad &y<1/2 \\ \frac{1}{Fr^2}\left(y \rho_{low} - \frac{\rho_{high} - \rho_{low}}{2}\right) + 1,\quad &y\geq 1/2 \end{cases},
\end{equation*}
and the ratio of specific heats is fixed at $\gamma = 5/3$. The Froude number is set as $Fr=1$, the Prandtl number remains $Pr=0.71$, and the Reynolds number is set to $Re=20{,}000$.

The instability is seeded by a small vertical velocity perturbation given as
\begin{equation}
  v = -0.025\sqrt{\frac{\gamma p}{\rho}}\cos(8\pi x).
\end{equation}
Finally, the $x-$direction boundaries are supplied with reflecting conditions, and the $y-$direction boundaries are held fixed at the initial density and pressure with zero velocity.

The source terms are constant and contain only quantities for which cell averages are already available, namely the density and $y-$momentum, and the averaged source $\langle\bm{S}\rangle$ is trivial to find for each cell. However, the implementation allows arbitrary user-provided source terms, so these gravitational sources are treated identically. As described in Section \ref{sec:pos_pres}, internal states can also be reconstructed within each cell on a tensor-product Gauss-Legendre set of nodes. The source term is evaluated over these states and subsequently integrated.

\begin{figure}
  \centering
  \includegraphics[width=0.5\linewidth]{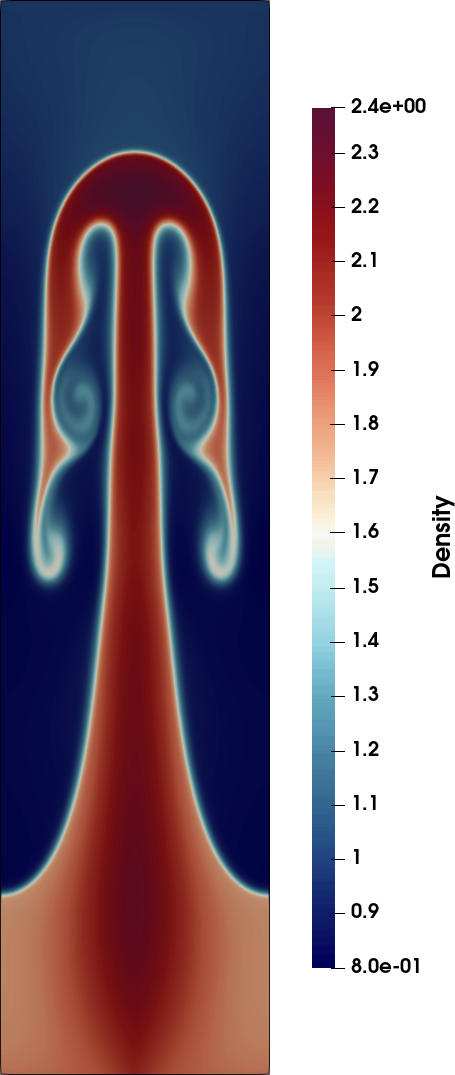}
  \caption{Shown is the density field for the viscous Rayleigh-Taylor instability at the final time of $t=2.5$ with Reynolds number $Re=20{,}000$ and Froude number $Fr=1$. These results were obtained with the radius $R=3$ scheme on a grid with spacing $\Delta=1/1024$, though they are insensitive to both grid resolution and to the stencil radius.}
  \label{fig:raytay}
\end{figure}

The problem is evolved to the final time of $t=2.5$ using the $[3S]_+^*$ time integrator with tolerances $atol=rtol=10^{-3}$ and a maximum CFL number of $1.25$. 
In \cref{fig:raytay}, we display the final density fields obtained by KFVM-WENO with radius $R=2$ and $R=3$ stencils on a grid spacing of $\Delta=1/1024$. Morphologically, with the explicit physical viscosity, a smooth leading cap is observed with no secondary Kelvin-Helmholtz type instabilities along its interface, and similarly there are no secondary instabilities along the central column. Secondary Kelvin-Helmholtz instabilities are visible on the inner region of the rising cap, and their structure is consistent over a range of grid resolutions. The two hooks (or arms) at the lowest part of the cap, as well as the position and the shape of the roll-up just above it, are converged and appear the same for all resolutions above $\Delta = 1/256$. The results are also consistent between the $R=2$ and $R=3$ schemes, with at most a $2\%$ difference in density in corresponding cells.

However, the above convergent solution behavior is in opposition to results obtained from inviscid solvers that exclude explicit viscosity but rely only on numerical dissipation because, in that case, there is no agreeable solution for the method to converge to. Indeed finer and finer scale structures will appear each time the grid is refined or the numerical dissipation is lowered by using a larger stencil (see \cite{shiRayTay:2003} for instance). Furthermore, methods with sufficiently low numerical dissipation are prone to breaking symmetry due to the non-associativity of floating point addition \citep{fleischmann2019numerical}. While fixes for this issue are available for dimension-by-dimension schemes, it remains unclear how one would avoid non-associativity errors in multidimensional reconstruction without drastically increasing the computational cost. Using a physical viscosity as done here avoids all of these problems by setting a single agreeable solution to converge to, and yields more scientifically meaningful results.

\section{Conclusion} \label{sec:conc}
This paper proposes a multidimensional adaptive order WENO finite volume method using kernel-based reconstruction. The use of kernel-based reconstruction allows great flexibility in the choice of stencils and substencils in multiple spatial dimensions. We showed in this paper that our non-polynomial, kernel-based design simplifies the implementation of high-order finite volume schemes in multiple dimensions by reconstructing all pointwise Riemann states along cell boundaries \textit{directly} from cell-average data, and eliminating the need to define boundary conditions for intermediate quantities as in modified dimension-by-dimension schemes.

Furthermore, alternative variables for reconstruction, dubbed linearized primitive variables, are proposed as a simplification over the use of characteristic variables. These are simpler to define and implement than characteristic variables. Crucially, these provide direction-independent information, which allows the same nonlinear weights within the WENO method to be used for the reconstruction of all Riemann states on all faces of a given cell. Alongside these new variables for reconstruction, we also proposed a straightforward adaptation of the KXRCF troubled cell indicator for use in finite volume schemes allowing WENO to be completely sidestepped in most cases. The calculation of nonlinear weights is the most expensive part of the whole scheme, so the use of these variables and the troubled cell indicator together provide a significant reduction in computational cost.

Finally, the proposed scheme is evaluated against a variety of stringent and illustrative benchmark problems. The method simultaneously demonstrates high-order nonlinear accuracy on smooth flows, robust behavior in the face of strong shocks, and minimal preference for grid-aligned phenomena over non-aligned phenomena.

\begin{acknowledgments}
This work was supported in part by the National Science Foundation under grants AST-1908834 and AST-2307684. We also acknowledge the use of the Lux supercomputer at UC Santa Cruz, funded by NSF MRI grant AST-1828315.
\end{acknowledgments}

\bibliography{new_merged}{}
\bibliographystyle{aasjournal}

\appendix

\section{Transformation matrices for linearized primitive variables}

Section \ref{sec:weno:vars:prim} presented the linearized primitive variables for reconstruction. For the sake of space the transformation matrices have been deferred to here. forward transformation matrix is
\begin{equation}
  \bm{\Phi} = \left.\frac{\partial\bm{V}}{\partial\bm{U}}\right|_{\langle\widetilde{\bm{U}}\rangle} = \begin{pmatrix} 1 & 0 & 0 & 0 & 0 \\ -\widetilde{u_1}/\widetilde{\rho} & 1/\widetilde{\rho} & 0 & 0 & 0 \\ -\widetilde{u_2}/\widetilde{\rho} & 0 & 1/\widetilde{\rho} & 0 & 0 \\ -\widetilde{u_3}/\widetilde{\rho} & 0 & 0 & 1/\widetilde{\rho} & 0 \\ (\gamma - 1)||\widetilde{u}||^2/2 & (1 - \gamma)\widetilde{u_1} & (1 - \gamma)\widetilde{u_2} & (1 - \gamma)\widetilde{u_3} & (\gamma - 1) \end{pmatrix},
  \label{eq:consToLPMat}
\end{equation}
where tildes indicate values obtained from the reference state $\langle\widetilde{\bm{U}}\rangle$. The velocities, $\widetilde{u_i} = \langle\widetilde{\rho u_i}\rangle/\langle\widetilde{\rho}\rangle$, are of course only second order accurate. The inverse transformation matrix is
\begin{equation}
  \bm{\Phi}^{-1} = \left.\frac{\partial\bm{U}}{\partial\bm{V}}\right|_{\bm{V}(\langle\widetilde{\bm{U}}\rangle)} = \begin{pmatrix} 1 & 0 & 0 & 0 & 0 \\ \widetilde{u_1} & \widetilde{\rho} & 0 & 0 & 0 \\ \widetilde{u_2} & 0 & \widetilde{\rho} & 0 & 0 \\ \widetilde{u_3} & 0 & 0 & \widetilde{\rho} & 0 \\ ||\widetilde{u}||^2/2 & \widetilde{u_1} & \widetilde{u_2} & \widetilde{u_3} & 1/(\gamma - 1) \end{pmatrix}.
  \label{eq:lpToConsMat}
\end{equation}
As mentioned in Section \ref{sec:weno:vars:prim}, both of these matrices contain a large number of zeros and the transformations can be applied more cheaply by ignoring these entries.

\end{document}